\documentclass[twocolumn]{autart}    

\usepackage{graphicx}          

\usepackage{amssymb}

\def\mve{\varepsilon}
\def\ms{\sigma}
\def\mtbb{\bm{\bar\theta}}
\def\mtb{\bar{\theta}}
\def\mt{\theta}
\def\Ex{{\rm I}\!{\rm E}}
\def\TT{^\top}
\def\ub{{\bf u}}
\def\mp{\partial}
\def\SN{{\mathcal N}}
\def\bm#1{\mbox{\boldmath $#1$}}
\def\eb{{\bf e}}
\def\Mb{{\bf M}}
\def\me{\epsilon}
\def\tb{{\bf t}}
\def\mthb{\bm{\hat\theta}}
\def\ra{\rightarrow}
\def\yb{{\bf y}}
\def\SD{{\mathcal D}}
\def\yh{\hat{y}}
\def\zb{{\bf z}}
\def\vb{{\bf v}}
\def\cb{{\bf c}}
\def\Cb{{\bf C}}
\def\mth{\hat{\theta}}
\def\1b{\mbox{\bf 1}}
\def\Ib{{\bf I}}
\def\rb{{\bf r}}
\def\SU{{\mathcal U}}
\def\Fx{{\rm I}\!{\rm F}}
\def\raas{\stackrel{\rm a.s.}{\ra}}
\def\rad{\stackrel{\rm d}{\ra}}
\def\ml{\lambda}
\def\SI{{\mathcal I}}
\def\ma{\alpha}
\def\mttb{\bm{\tilde\theta}}
\def\SR{{\mathcal R}}
\def\SX{{\mathcal X}}
\def\tr{\mbox{\rm trace}}
\def\Qb{{\bf Q}}
\def\SM{{\mathcal M}}
\def\var{\mbox{\rm var}}
\def\Hb{{\bf H}}
\newtheorem{theo}{\bf Theorem}
\def\Var{\mbox{\rm Var}}
\def\0b{\mbox{\bf 0}}
\def\SP{{\mathcal P}}
\def\SO{{\mathcal O}}
\def\SF{{\mathcal F}}
\def\Ob{{\bf O}}

\begin{document}

\begin{frontmatter}
\runtitle{Optimal design and control}  

\title{Optimal experimental design\\ and some related control problems\thanksref{footnoteinfo}} 

\thanks[footnoteinfo]{This paper was not presented at any IFAC
meeting. Corresponding author L.~Pronzato. Tel. +33 (0)4 92942708.
Fax +33 (0)4 92942896.}


\author{Luc Pronzato}\ead{pronzato@i3s.unice.fr},
\address{Laboratoire I3S, CNRS-Universit\'e de Nice-Sophia Antipolis, France}


\begin{abstract}                          
This paper traces the strong relations between experimental design
and control, such as the use of optimal inputs to obtain precise
parameter estimation in dynamical systems and the introduction of
suitably designed perturbations in adaptive control. The
mathematical background of optimal experimental design is briefly
presented, and the role of experimental design in the asymptotic
properties of estimators is emphasized. Although most of the paper
concerns parametric models, some results are also presented for
statistical learning and prediction with nonparametric models.
\end{abstract}

\end{frontmatter}

\section{Introduction}

The design of experiments (DOE) is a well developed methodology in
statistics, to which several books have been dedicated, see e.g.\
\cite{Fedorov72}, \cite{Silvey80}, \cite{Pazman86},
\cite{AtkinsonD92}, \cite{Pukelsheim93}, \cite{FedorovH97}. See
also the series of proceedings of the Model-Oriented Design and
Analysis workshops (Springer Verlag 1987; Physica Verlag, 1990,
1992, 1995, 1998, 2001, 2004). Its application to the construction
of persistently exciting inputs for dynamical systems is well
known in control theory (see Chapter 6 of \cite{GoodwinP77},
Chapter 14 of \cite{Ljung87}, Chapter 6 of \cite{WPl97}, the book
\cite{Zarrop79} and the recent surveys \cite{Gevers2005},
\cite{Hjalmarsson2005}). A first objective of this paper is to
briefly present the mathematical background of the methodology and
make it accessible to a wider audience. DOE, which can can be
apprehended as a technique for extracting the most useful
information from data to be collected, is thus {\em a central (and
sometimes hidden) methodology in every occasion where unknown
quantities must be estimated and the choice of a method for this
estimation is open}. DOE may therefore serve different purposes
and happens to be a suitable vehicle for establishing links
between problems like optimization, estimation, prediction and
control. Hence, a second objective of the paper is to exhibit
links and similarities between seemingly different issues (for
instance, we shall see that parameter estimation and prediction of
a model response are essentially equivalent problems for
parametric models and that the construction of an optimal method
for global optimization can be casted as a stochastic control
problem). At the same time, attention will be drawn to fundamental
differences that exist between seemingly similar problems (in
particular, evidence will be given of the gap between using
parametric or nonparametric models for prediction). A third
objective is to point out and explain some inherent difficulties
in estimation problems when combined with optimization or control
(hence we shall see why adaptive control is an intrinsically
difficult subject), indicate some tentative remedies and suggest
possible developments.

Mentioning these three objectives should not shroud the main
message of the paper, which consists in {\em pointing out
prospective research directions for experimental design in
relation with control}, in short: classical DOE relies on the
assumption of persistence of excitation but many issues remain
open in other situations; DOE should be driven by the final
purpose of the identification (the intended model application of
\cite{GeversL86}) and this should be reflected in the construction
of design criteria; DOE should face the new challenges raised by
nonparametric models and robust control; algorithms and practical
methods for DOE in non-standard situations are still missing. The
program is rather ambitious, and this survey does not pretend to
be exhaustive (for instance, only the case of scalar observations
is considered; Bayesian techniques are only slightly touched;
measurement errors are assumed to be independent, although
correlated errors would deserve a special treatment; distributed
parameter systems are not considered; nonparametric modelling is
briefly considered and for static systems only, etc.). However,
references are indicated where a detailed enough presentation is
lacking. None of the results presented is really new, but their
collection in a single document probably is, and will hopefully be
useful to the reader.

Section \ref{S:examples} presents different types of application
of optimal experimental design, partly through examples, and
serves as an introduction to the topic. In particular, the fourth
application concerns optimization and forms a preliminary
illustration of the link between sequential design and adaptive
control. Section \ref{nonparametric} concerns statistical learning
and nonparametric modelling, where DOE is still at an early stage
of development. The rest of the paper mainly deals with parametric
models, for which parameter uncertainty is suitably characterized
through information matrices, due to the asymptotic normality of
parameter estimators and the Cram\'er-Rao bound. This is
considered in Section \ref{S:information-matrix} for regression
models. Section \ref{S:design-for-estimation} presents the
mathematical background of optimal experimental design for
parameter estimation. The case of dynamical models is considered
in Section \ref{S:control-for-design}, where the input is designed
to yield the most accurate estimation of the model parameters,
while possibly taking a robust-control objective into account.
Section \ref{S:adaptive-control} concerns adaptive control: the
ultimate objective is process control, but the construction of the
controller requires the estimation of the model parameters. The
difficulties are illustrated through a series of simple examples.
Optimal DOE yields input sequences that are optimally (and
persistently) exciting. At the same time, by focussing attention
on parameter estimation, it reveals the intrinsic difficulties of
adaptive control through the links between dual (active) control
and sequential design. General sequential design (for static
systems) is briefly considered in Section
\ref{S:sequential-design}. Finally, Section \ref{S:further-issues}
suggests further developments and research directions in DOE,
concerning in particular active learning and nonlinear feedback
control. Here also the presentation is mainly through examples.

\section{Examples of applications of DOE} \label{S:examples}

Although the paper is mainly dedicated to parameter estimation
issues, DOE may have quite different objectives (and it is indeed
one of the purposes of the paper to use DOE to exhibit links
relating these objectives). They are illustrated through examples
which also serve to progressively introduce the notations. The
first one concerns an extremely simple parameter estimation
problem where the benefit of a suitably designed experiment is
spectacular.

\subsection{A weighing problem}\label{S:weighing}

Suppose we wish to determine the weights of eight objects with a
chemical balance. The result $y$ of a weighing (the observation)
corresponds to the mass on the left pan of the balance minus the
mass on the right pan plus some measurement error $\mve$. The
errors associated with a series of measurements are assumed to be
independently identically distributed (i.i.d.) with the normal
distribution $\SN(0,\ms^2)$. The objects have weights $\mtb_i$,
$i=1,\ldots,8$. Each weighing is characterized by a 8-dimensional
vector $\ub$ with components $\{\ub\}_i$ equal to $1$, $-1$ or $0$
depending whether object $i$ is on the left pan, the right pan or
is absent from the weighing, and the associated observation is
$y=\ub\TT\mtbb + \mve$. We thus have a linear model (in the
statistical sense: the response is linear in the parameter vector
${\bm\mt}$), and the Least-Squares (LS) estimator $\bm{\hat\mt}^N$
associated with $N$ observations $y_k$ characterized by
experimental conditions (design points\footnote{Although design
points and experimental variables are usually denoted by the
letter $x$ in the statistical literature, we shall use the letter
$u$ due to the attention given here to control problems. In this
weighing example, $\ub_k$ denotes the decisions made concerning
the $k$-th observation, which already reveals the contiguity
between experimental design and control.}) $\ub_k$,
$k=1,\ldots,N$, is
\begin{eqnarray}
&& \bm{\hat\mt}^N = \arg\min_{\bm\mt} \sum_{k=1}^N
[y_k-\ub_k\TT{\bm\mt}]^2 = \Mb_N^{-1} \sum_{k=1}^N
        y_k\ub_k \,, \label{LSlinear} \\
&& \mbox{with }
    \Mb_N = \sum_{k=1}^N \ub_k \ub_k\TT \,. \label{Mlinear}
\end{eqnarray}
We consider two weighing methods.
In method $a$ the eight objets are weighed successively: the
vectors $\ub_i$ for the eight observations coincide with the basis
vectors $\eb_i$ of $\mathbb{R}^8$ and the observations are
$y_i=\mtb_i+\mve_i$, $i=1,\ldots,8$. The estimated weights are
simply given by the observations, that is, $\hat
\mt_i=y_i\sim\SN(\mtb_i,\ms^2)$. Method $b$ is slightly more
sophisticated. Eight measurements are performed, each time using a
different configuration of the objets on the two pans so that the
vectors $\ub_i$ form a $8\times 8$ Hadamard matrix
($|\{\ub_i\}_j|=1$ $\forall i,j$ and $\ub_i\TT\ub_j=0$ $\forall
i\neq j$, $i,j=1,\ldots,8$). The estimates then satisfy $\hat
\mt_i\sim\SN(\mtb_i,\ms^2/8)$ with 8 observations only. To obtain
the same precision with method $a$, one would need to perform
eight independent repetitions of the experiment, requiring 64
observations in total\footnote{Note that we implicitly assumed
that the range of the instrument allows to weigh all objects
simultaneously in method $b$. Also note that the gain would be
smaller when using method $b$ if the variance of the measurement
errors increased with the total weight on the balance.}.

In a linear model of this type, the LS estimator (\ref{LSlinear})
is unbiased: $\Ex_{\mt}\{\bm{\hat\mt}^N \} - {\bm\mt} = 0$, where
$\Ex_{\mt}\{\cdot\}$ denotes the mathematical expectation
conditionally to ${\bm\mt}$ being the true vector of unknown
parameters. Its covariance matrix is $\Ex_{\mt}\{
(\bm{\hat\mt}^N-{\bm\mt})(\bm{\hat\mt}^N-{\bm\mt})\TT \} = \ms^2
\Mb_N^{-1}$ with $\Mb_N$ given by (\ref{Mlinear}) (note that it
does not depend on $\mt$). Choosing an experiment that provides a
precise estimation of the parameters thus amounts to choosing $N$
vectors $\ub_k$ such that ($\Mb_N$ is non singular and)
``$\Mb_N^{-1}$ is as small as possible'', in the sense that a
scalar function of $\Mb_N^{-1}$ is minimized (or a scalar function
of $\Mb_N$ is maximized), see Section
\ref{S:design-for-estimation}. In the weighing problem above the
optimization problem is combinatorial since $\{\ub_k\}_i \in
\{-1,0,1\}$. In the design of method $b$ the vectors $\ub_k$
optimize most ``reasonable'' criteria $\Phi(\Mb_N)$, see, e.g.,
\cite{Chen80}, \cite{Schwabe87b}. This case will not be considered
in the rest of the paper but corresponds to a topic that has a
long and rich history (it originated in agriculture through the
pioneering work of Fisher, see \cite{Fisher25}).

\subsection{An example of parameter estimation in a dynamical
model}\label{S:pharmaco}

\begin{figure}[h]
\begin{center}
\includegraphics[scale=0.35]{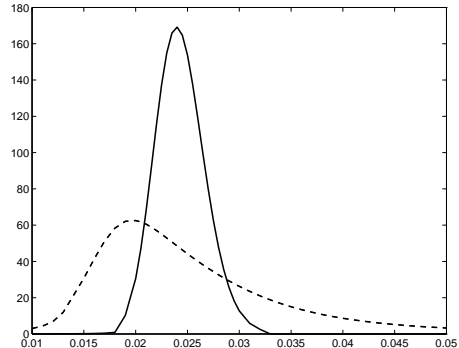}
\caption{Approximate marginal densities for the LS estimator
$\hat{K}_{EL}$ (dashed line for the conventional design, solid
line for the optimal one); the true value is $\bar
K_{EL}=0.0242\mbox{ min}^{-1}$.} \label{F:q1}
\end{center}
\end{figure}

The example is taken from \cite{D'Argenio81} and concerns a
so-called compartment model, widely used in pharmacokinetics. A
drug $x$ is injected in blood (intravenous infusion) with an input
profile $u(t)$, the drug moves from the central compartment $C$
(blood) to the peripheral compartment $P$, where the respective
quantities of drugs at time $t$ are denoted $x_C(t)$ and $x_P(t)$.
These obey the following differential equations:
$$
\left\{ \begin{array}{rcllcl}
\frac{d x_C(t)}{dt} & = & (- K_{EL} - & K_{CP}) x_C(t) & + & K_{PC} x_P(t) + u(t) \\
\frac{d x_P(t)}{dt} & = & & K_{CP} x_C(t) & - & K_{PC} x_P(t) \\
\end{array} \right.
$$
where $K_{CP}$, $K_{PC}$ and $K_{EL}$ are unknown parameters. One
observes the drug concentration in blood, that is, $y(t)=x_C(t)/V
+ \mve(t)$ at time $t$, where the errors $\mve(t_i)$ corresponding
to different observations are assumed to be i.i.d.\ $\SN(0,\ms^2)$
and where $V$ denotes the (unknown) volume of the central
compartment. There are thus four unknown parameters to be
estimated, which we denote $\bm\mt=(K_{CP},K_{PC},K_{EL},V)$. The
profile of the input $u(t)$ is imposed: it consists of a 1 min
loading infusion of 75 mg/min followed by a continuous maintenance
infusion of 1.45 mg/min. The experimental variables correspond to
the sampling times $t_i$, $1\leq t_i\leq 720 \mbox{ min}$ (the
time instants at which the observations ---~blood samples~--- are
taken). Suppose that the true parameters take the values $\mtbb =
(0.066\, \mbox{min}^{-1},\ 0.038\, \mbox{min}^{-1}, \ 0.0242\,
\mbox{min}^{-1}, \ 30\, \mbox{l)}$. Two different experimental
designs are considered. The first one, called ``conventional'', is
given by $\tb = (5, 10, 30, 60, 120, 180, 360, 720) \mbox{ (in
min)}$; the ``optimal'' one ($D$-optimal for $\mtbb$, see Section
\ref{S:criteria}) is $\tb^* = (1, 1, 10, 10, 74, 74, 720, 720)
\mbox{ (in min)}$. (Note that both designs contain 8 observations
and that $\tb^*$ comprises repetitions of observations at the same
time ---~which means that it is implicitly assumed that the
collection of several simultaneous independent measurements is
possible.) Figure \ref{F:q1} presents the (approximate) marginal
density for the LS estimator of $K_{EL}$, see \cite{PazmanPa96},
\cite{PP01a}, when $\ms=0.2 \mu\mbox{g/ml}$. Similar pictures are
obtained for the other parameters.

Clearly, the ``optimal'' design $\tb^*$ yields a much more precise
estimation of $\bm\mt$ than the conventional one, although both
involve the same number of observations. On the other hand, with
$4=\dim(\bm\mt)$ sampling times only, $\tb^*$ does not permit to
test the validity of the model. DOE for model discrimination,
which we consider next, is especially indicated for situations
where one hesitates between several structures.

\subsection{Discrimination between model structures}
\label{S:discrimination}


Design for discrimination between model structures will not be
detailed in the paper, only the basic principle of a simple method
is indicated below and one can refer to \cite{BoxH67} and the
survey papers  \cite{AtkinsonC74}, \cite{Hill78} for other
approaches. When there are two model structures
$\eta^{(1)}(\bm\mt_1,\ub)$ and $\eta^{(2)}(\bm\mt_2,\ub)$ and the
errors are i.i.d., a simple
sequential procedure is as follows, see \cite{AtkinsonF75a}:\\
$\bullet$ after the observation of $y(\ub_1),\ldots,y(\ub_k)$
estimate $\mthb_1^k$ and $\mthb_2^k$ for both
    models;\\
$\bullet$ place next point $\ub_{k+1}$ where
$[\eta^{(1)}(\mthb_1^k,\ub)-\eta^{(2)}(\mthb_2^k,\ub)]^2$ is
    maximum;\\
$\bullet$ $k \ra k+1$, repeat. \\
When there are more than two structures in competition, one should
estimate $\mthb_i^k$ for all of them and place the next point
using the two models with the best and second best fitting, see
\cite{AtkinsonF75b}. The idea is to place the design point where
the predictions of the competitors differ much, so that when one
of the structures is correct (which is the underlying assumption),
next observation should be close to the prediction of that model
and should thus give evidence that the other structures are wrong.
Similar ideas can be used to design input sequences for detecting
changes in the behavior of dynamical systems, see the book
\cite{Kerestecioglu93}.

\subsection{Optimization of a model
response}\label{S:optimization}

Suppose that one wishes to maximize a function $\eta(\mtbb,\ub)$
with respect to $\ub\in\mathbb{R}^d$, with $\mtbb\in\mathbb{R}^p$
a vector of unknown parameters. When a value $\ub_i$ is proposed,
the function is observed through
$y_i=y(\ub_i)=\eta(\mtbb,\ub_i)+\mve_i$ with $\mve_i$ a
measurement error. Since the problem is to determine
$\ub^*=\ub^*(\mtbb)=\arg\max_{\ub} \eta(\mtbb,\ub)$, it seems
natural to first estimate $\mthb=\mthb[\yb]$ from a vector of
observations $\yb=[y_1,\ldots,y_N]\TT$ and then predict the
optimum by $\ub^*(\mthb)$. The question is then which values to
use for the $\ub_i$'s for estimating $\mthb$, that is, which
criterion to optimize for designing the experiment? It could be
(i) based on the precision of $\mthb$, or (ii) based on the
precision of $\ub^*(\mthb)$, or, preferably, (iii) {\em oriented
towards the final objective} and based on the cost
$C(\mthb|\mtbb)$ of using $\mthb$ when the true value of $\bm\mt$
is $\mtbb$. A possible choice is $C(\mthb|\mtbb) =
\eta[\mtbb,\ub^*(\mtbb)] - \eta[\mtbb,\ub^*(\mthb)] \geq 0$, which
leads to a design that minimizes the Bayesian risk $R = \Ex
\{C(\mthb[\yb]|\mtbb)\}$, where the expectation is with respect to
$\yb$ and $\mtbb$ for which a prior distribution $\pi(\cdot)$ is
assumed, see, e.g., \cite{PWa93} (see also \cite{ChalonerV95} and
the book \cite{Pilz83} for a review of Bayesian DOE).

The approaches (i-iii) above are standard in experimental design:
optimization is performed in two steps, first some design points
$\ub_i$'s are selected for estimation, second $\mthb$ is estimated
and used to construct $\ub^*(\mthb)$. However, in general each
response $\eta(\mtbb,\ub_i)$ is far from the maximum
$\eta[\mtbb,\ub^*(\mtbb)]$ (since the explicit objective of the
design is estimation, not maximization) while in some situations
it is required to have $\eta(\mtbb,\ub_i)$ as large as possible
for every $i$, that is, $\ub_i$ close to $\ub^*(\mtbb)$, which is
unknown. A sequential approach is then natural: try $\ub_i$,
observe $y_i$, estimate $\mthb^i=\mthb(\yb_1^i)$, suggest
$\ub_{i+1}$ and so on\ldots (Notice that this involves a feedback
of information in the sequence of design points ---~the control
sequence~--- and thus induces a dynamical aspect although the
initial problem is purely static.) Each $\ub_i$ has two
objectives: help to estimate $\bm\mt$, try to maximize
$\eta(\mtbb,\ub)$. The design problem thus corresponds to a {\em
dual control} problem, to be considered in Section
\ref{S:finiteN}. When no parametric form is known for the function
to be maximized, it is classical to resort to suboptimal methods
such as the Kiefer-Wolfowitz scheme \cite{KieferW52}, or the
response surface methodology which involves linear and quadratic
approximations, see, e.g., \cite{BoxW51}. Optimization with a
nonparametric model will be considered in Section
\ref{S:further-issues}, combining statistical learning with global
optimization.

\section{Statistical learning, nonparametric models}\label{nonparametric}

One can refer to the books \cite{Vapnik98}, \cite{Vapnik2000},
\cite{HastieTF2001} and the surveys \cite{CuckerS2001},
\cite{Bartlett2003} for a detailed exposition of statistical
learning. Based on so-called ``training data''
$\SD=\{[\ub_1,y(\ub_1)],\ldots,$ $[\ub_N,y(\ub_N)]\}$ we wish to
predict the response $y(\ub)$ of a process at some unsampled input
$\ub$ using Nadaraya-Watson regression \cite{Nadaraya64},
\cite{Watson64}, Radial Basis Functions (RBF), Support Vector
Machine (SVM) regression or Kriging (Gaussian process). All these
approaches can be casted in the class of kernel methods, see
\cite{Vazquez05}, \cite{VazquezW03} and \cite{Schaback2003} for a
more precise formulation, and we only consider the last one,
Kriging, due to its wide flexibility and easy interpretability.
The associated DOE problem is considered in Section
\ref{S:DOEnonparametric}. We denote $\yh_\SD(\ub)$ the prediction
at $\ub$ and $\yb=[y(\ub_1),\ldots,y(\ub_N)]\TT$.

\subsection{Gaussian process and Kriging}\label{S:Kriging}

The method originated in geostatistics, see \cite{Krige51},
\cite{Matheron63}, and has a long history. When the modelling
errors concern a transfer function observed in the Nyquist plane,
the approach possesses strong similarities with the so-called
``stochastic embedding'' technique, see, e.g., \cite{GoodwinS89}
and the survey paper \cite{NinnessG95}. The observations are
modelled as $y(\ub_k)= \mtb_0 + P(\ub_k,\omega) + \mve_k$, where
$P(\ub,\omega)$ denotes a second-order stationary zero-mean random
process with covariance
$\Ex\{P(\ub,\omega)P(\zb,\omega)\}=K(\ub,\zb)=\ms^2_P C(\ub-\zb)$
and the $\mve_k$'s are i.i.d., with zero mean and variance
$\ms^2$. The best linear unbiased predictor at $\ub$ is
$\yh_{\SD}(\ub)=\vb\TT(\ub) \yb$, where $\vb(\ub)$ minimizes
$\Ex\{ (\vb\TT\yb - [\mtb_0 + P(\ub,\omega)])^2\}$ with the
constraint $\Ex\{\vb\TT\yb\}=\mtb_0 \sum_{i=1}^N
v_i=\Ex\{y(\ub)\}=\mtb_0$, that is, $\sum_{i=1}^N v_i=1$. This
optimization problem is solvable explicitly, which gives
\begin{equation}\label{yh-kriging}
    \yh_{\SD}(\ub)=\vb\TT(\ub)\yb= \mth^0 + \cb\TT(\ub)
\Cb_y^{-1}(\yb-\mth^0 \1b)
\end{equation}
where $\Cb_y = \ms^2 \Ib_N + \ms^2_P \Cb_P$ with $\Ib_N$ the
$N$-dimensional identity matrix and $\Cb_P$ the $N\times N$ matrix
defined by $[\Cb_P]_{i,j}=C(\ub_i-\ub_j)$, $\1b$ is the
$N$-dimensional vector with components 1, $\cb(\ub)=\ms^2_P
[C(\ub-\ub_1),\ldots,C(\ub-\ub_N)]\TT$ and $\mth_0= (\1b\TT
\Cb_y^{-1}\yb)/(\1b\TT \Cb_y^{-1} \1b)$ (a weighted LS estimator
of $\mt_0$). Note that the prediction takes the form
$\yh_{\SD}(\ub)=\mth^0 + \sum_{k=1}^N a_k K(\ub,\ub_k)$, i.e., a
linear combination of kernel values. The Mean-Squared-Error (MSE)
of the prediction $\yh_\SD(\ub)$ at $\ub$ is given by
\begin{equation}\label{MSE-kriging}
 \rho_\SD^2(\ub)=\ms^2_P- \left[
\begin{array}{cc}
  \cb\TT(\ub) & 1 \\
\end{array} \right]
\left[ \begin{array}{cc}
  \Cb_y & \1b \\
  \1b\TT & 0 \\
\end{array} \right]^{-1}
\left[ \begin{array}{c}
  \cb(\ub) \\
  1 \\
\end{array} \right]
\end{equation}
and, if $\ms^2=0$ (i.e., there are no measurement errors
$\mve_k$), $\yh_{\SD}(\ub_i)=y(\ub_i)$ and $\rho_\SD^2(\ub_i)=0$
for any $i$. The predictor $\yh_{\SD}(\ub)$ is then a perfect
interpolator. This method thus makes statistical inference
possible even for purely deterministic systems, the model
uncertainty being represented by the trajectory of a random
process. Since the publication \cite{SacksWMW89} it has been
successfully applied in many domains of engineering where
simulations (computer codes) replace real physical experiments
(and measurement errors are thus absent), see, e.g.,
\cite{SantnerWN2003}.

%

If the characteristics of the process $P(\ub,\omega)$ and errors
$\mve_k$ belong to a parametric family, the unknown parameters
that are involved can be estimated. For instance, for a Gaussian
process with $C(\zb)$ parameterized as $C(\zb)=C(\bm\beta,\zb)$
and for normal errors $\mve_k$, the parameters $\bm\beta$,
$\ms^2_P$ and $\ms^2$ can be estimated by Maximum Likelihood; see
the book \cite{Stein99}, in particular for recommendations
concerning the choice of the covariance function $C(\zb)$. See
also the survey \cite{Mardia90} and the papers \cite{Ying93},
\cite{vanderVaart96} concerning the asymptotic properties of the
estimator. The method can be extended in several directions: the
constant terms $\mt_0$ can be replaced by a linear model
$\rb\TT(\ub)\bm\mt$ (this is called universal Kriging, or {\em
intrinsic Kriging} when generalized covariances are used, which is
then equivalent to {\em splines}, see \cite{VazquezWF2005}), a
prior distribution can be set on $\bm\mt$ (Bayesian Kriging, see
\cite{CurrinMMY91}), the derivative (gradient) of the response
$y(\ub)$ can also be predicted from observations $y(\ub_k)$, see
\cite{Vazquez05}, or observations of the derivatives can be used
to improve the prediction of the response, see \cite{MorrisMY93},
\cite{MardiaKGL96}, \cite{LearyBK2004}. Nonparametric modelling
can be used in optimization, and an application of Kriging to
global optimization is presented in Section
\ref{S:further-issues}.

\subsection{DOE for nonparametric
models}\label{S:DOEnonparametric}

The approaches can be classified among those that are model-free
(of the space-filling type) and those that use a model.

\subsubsection{Model-free design (space filling)}

For $\SU$ the design set (the admissible set for $\ub$), we call
$\SS\subset\SU$ the finite set of chosen design points or sites
$\ub_k$ where the observations are made, $k=1,\ldots,N$. {\em
Maximin-distance design} \cite{JohnsonMY90} chooses sites $\SS$
that maximize the minimum distance between points of $\SS$, i.e.\
$\min_{\ub\neq \ub'\in\SS^2} d(\ub,\ub')$. The chosen sites
$\ub_k$ are thus maximally spread in $\SU$ (in particular, some
points are set on the boundary of $\SU$). When $\SU$ is a discrete
set, {\em minimax-distance design} \cite{JohnsonMY90} chooses
sites that minimize the maximum distance between a point in $\SU$
and $\SS$, i.e.\ $\max_{\zb\in\SU} d(\zb,\SS)  = \max_{\zb\in\SU}
\min_{\ub\in\SS} d(\zb,\ub)$. In order to ensure good projection
properties in all directions (for each component of the
$\ub_k$'s), it is recommended to work in the class of {\em latine
hypercube designs}, see \cite{MorrisM95} (when $\SU$ is scaled to
$[0,1]^d$, for every $i=1,\ldots,d$ the components $\{\ub_k\}_i$,
$k=1,\ldots,N$, then take all the values
$0,1/(N-1),2/(N-1),\ldots,1$).

\subsubsection{Model-based design}\label{S:Model-based}

In order to relate the choice of the design to the quality of the
prediction $\yh_D(\ub)$, a first step is to characterize the
uncertainty on $\yh_D(\ub)$. This raises difficult issues in
nonparametric modelling, in particular due to the difficulty of
deriving a global measure expressing the speed of decrease of the
MSE of the prediction as $N$, the number of observations,
increases (we shall see in Section \ref{S:active-parametric} that
the situation is opposite in the parametric case). A reason is
that the effect of the addition of a new observation is local:
when we observe at $\ub$, the MSE of the prediction at $\zb$
decreases for $\zb$ close to $\ub$ (for instance, for Kriging
without measurement errors $\rho_\SD(\ub)$ becomes zero), but is
weakly modified for $\zb$ far from $\ub$. Hence, DOE is often
ignored in the statistical learning literature\footnote{There
exists a literature on {\em active learning}, which aims at
selecting training data using techniques from DOE. However, it
seems that when explicit reference to DOE is made, the attention
is restricted to learning with a parametric model, see in
particular \cite{Cohn94}, \cite{CohnGJ96}. In that case, the
underlying assumption that the data are generated by a process
whose structure coincides with that of the model is often hardly
tenable, especially for a behavioral model e.g.\ of the
neural-network type, see Section \ref{S:active-parametric} for a
discussion.}, where the set of training data $\SD$ is generally
assumed to be a collection of i.i.d.\ pairs $[\ub_k,y(\ub_k)]$,
see, e.g., \cite{CuckerS2001}, \cite{Bartlett2003}. The local
influence just mentioned has the consequence that an optimal
design should (asymptotically) tend to observe everywhere in
$\SU$, and distribute the points $\ub_k$ with a density (i.e.\
according to a probability measure absolutely continuous with
respect to the Lebesgue measure on $\SU$ ---~again, we shall see
that the situation is opposite for the parametric case). Few
results exist on that difficult topic, see e.g.\ \cite{ChengHT98}:
for $u$ scalar, observations $y(u_k)=f(u_k)+\mve_k$ with i.i.d.
errors $\mve_k$, and a prediction of the Nadaraya-Watson type
(\cite{Nadaraya64}, \cite{Watson64}), a sequential algorithm is
constructed that is asymptotically optimal (it tends to distribute
the points $u_k$ with a density proportional to $|f''(u)|^{2/9}$).
See also \cite{Muller84}, \cite{Faraway90} for related results.
The uniform distribution may turn out to be optimal when
considering minimax optimality over a class of functions, see
\cite{BiedermannD2001}.

When Kriging is used for prediction, the MSE is given by
(\ref{MSE-kriging}) and $\SS$ can be chosen for instance by
minimizing the maximum MSE
 $\max_{u\in\SU} \rho_\SD^2(\ub)$ (which is related to minimax-distance design,
 see \cite{JohnsonMY90}) or by minimizing the integrated MSE $\int_\SU \rho_\SD^2(\ub) \,
 \pi(d\ub)$, with $\pi(\cdot)$ some probability density for $\ub$, see
\cite{SacksS88}. {\em Maximum entropy sampling} \cite{ShewryW87}
provides an elegant alternative design method, usually requiring
easier computations. It can be related to maximin-distance design,
see \cite{JohnsonMY90}.

Notice finally that in general the parameters $\bm\beta$,
$\ms^2_P$ and $\ms^2$ in the covariance matrix $\Cb_y$ used in
Kriging are estimated from data, so that the precision of their
estimation influences the precision of the prediction. This seems
to have received very little attention, although designs for
prediction (space filling for instance) are clearly not
appropriate for the precise estimation of these parameters, see
\cite{ZhuZ2006}.

\section{Parametric models and information matrices} \label{S:information-matrix}

Throughout this section we consider regression models with
observations
\begin{equation}\label{regression}
    y(\ub_k)=\eta(\mtbb,\ub_k)+\mve_k\,, \ \mtbb\in\Theta\,, \
\ub_k\in\SU \,,
\end{equation}
where the errors $\mve_k$ are independent with zero mean and
variance $\Ex_{\ub_k}(\mve_k^2)=\ms^2(\ub_k)$, $k=1,2,\ldots$
(with $0<a\leq \ms^2(\ub)\leq b<\infty$). The function
$\eta(\bm\mt,\ub_k)$ is known, possibly nonlinear in $\bm\mt$, and
$\mtbb$, the true value of the model parameters, is unknown. The
asymptotic behavior of the LS estimator, in relation with the
design, is recalled in the next section (precise proofs are
generally rather technical, and we give conditions on the design
that facilitate their construction). Maximum-Likelihood estimation
and estimating functions are considered next. The extension to
dynamical systems requires more technical developments beyond the
scope of this paper. One can refer e.g.\ to \cite{GoodwinP77},
\cite{Ljung87}, \cite{Caines88} \cite{SoderstromS89} for a
detailed presentation, including data-recursive estimation
methods. Also, one can refer to the monograph \cite{Ucinski2005}
for the identification of systems with distributed parameters and
e.g.\ to \cite{KubruslyM85}, \cite{Rafajlowicz83},
\cite{Rafajlowicz86b}
for optimal input design for such systems.

\subsection{Weighted LS estimation}\label{S:WLS}

Consider the weighted LS (WLS) estimator
$$
    \mthb_{WLS}^N = \arg\min_{\bm\mt} (1/N)\sum_{k=1}^N w(\ub_k) [y(\ub_k) -
\eta(\bm\mt,\ub_k)]^2
$$
with $w(\cdot)$ a known function, bounded on $\SU$. To investigate
the asymptotic properties of $\mthb_{WLS}^N$ for $N\ra\infty$ we
need to specify how the design points $\ub_k$'s are generated. In
that sense, {\em the asymptotic properties of the estimator are
strongly related to the design}. The early and now classical
reference \cite{Jennrich69} makes assumptions on the finite tail
products of the regression function $\eta$ and its derivatives,
but the results are more easily obtained at least in
two cases: \\
$\bullet$ (i) $(\ub_k)$ forms a sequence of i.i.d.\ random
variables (vectors), distributed with a probability measure $\xi$
(which we
call {\em random design}); \\
$\bullet$ (ii) The empirical measure $\xi_N$ with distribution
function $\Fx_{\xi _N}(\ub)= \sum_{i=1,\,\,\ub_i < \ub}^N (1/N)$
(where the inequality $\ub_i<\ub$ is componentwise) converges
strongly (in variation, see \cite{Shiryaev96}, p.\ 360) to a
discrete probability measure $\xi$ on $\SU$, with finite support
$\SS_\xi =\left\{ \ub\in\SU: \xi (\{ \ub\})>0\right\}$, that is,
$\lim_{N\ra\infty} \xi_N(\{\ub\})= \xi(\{\ub\})$ for any
$\ub\in\SU$.

Note that in case (i) the pairs $(\mve_k,\ub_k)$ are i.i.d. and in
case (ii) there exist a finite number of {\em support points}
$\ub^i$ that receive positive weights $\xi(\ub^i)>0$, so that, as
$N$ increases, the observations at those $\ub^i$'s are necessarily
repeated. In both cases the asymptotic distribution of the
estimator is characterized by $\xi$.

The strong consistency of $\mthb_{WLS}^N$, i.e., $\mthb_{WLS}^N
\raas \mtbb$, $N\ra\infty$, can easily be proved for designs
satisfying (i) or (ii) under continuity and boundedness
assumptions on $\eta(\bm\mt,\ub)$ when the estimability condition
$[\int_\SU w(\ub)[\eta(\bm\mt,\ub)-\eta(\bm\mt',\ub)]^2\,
\xi(d\ub)=0 \Leftrightarrow \theta'=\theta]$ is satisfied.
Supposing, moreover, that $\eta(\bm\mt,\ub)$ is two times
continuously differentiable in $\bm\mt$ and that the matrix
$$ \Mb_1(\xi,\mtbb)= \int_\SU w(\ub)\,
\frac{\mp\eta(\bm\mt,\ub)}{\mp\bm\mt}_{|\mtbb}
\frac{\mp\eta(\bm\mt,\ub)}{\mp\bm\mt\TT}_{|\mtbb} \, \xi(d\ub)
$$
has full rank, an application of the Central Limit Theorem to a
Taylor series development of $\nabla_\mt J_N(\bm\mt)$, the
gradient of the WLS criterion, around $\mthb_{WLS}^N$ gives
\begin{equation}\label{ANWLS}
    \sqrt{N} (\mthb^N_{WLS}-\mtbb) \rad z \sim\SN(0,\Cb(w,\xi,\mtbb))\,,
\ N\ra\infty \,,
\end{equation}
where
$\Cb(w,\xi,\bm\mt)=\Mb_1^{-1}(\xi,\bm\mt)\Mb_2(\xi,\bm\mt)\Mb_1^{-1}(\xi,\bm\mt)$
with
$$
\Mb_2(\xi,\bm\mt)= \int_\SU w^2(\ub)\ms^2(\ub)\,
\frac{\mp\eta(\bm\mt,\ub)}{\mp\bm\mt}
\frac{\mp\eta(\bm\mt,\ub)}{\mp\bm\mt\TT}\, \xi(d\ub) \,.
$$
One may notice that $\Cb(w,\xi,\mtbb)-\Mb^{-1}(\xi,\mtbb)$ is
non-negative definite for any weighting function $w(\cdot)$, where
$\Mb(\xi,\bm\mt)$ denotes the matrix
\begin{equation}\label{MinfoWLS}
    \Mb(\xi,\bm\mt)= \int_\SU \ms^{-2}(\ub) \,
\frac{\mp\eta(\bm\mt,\ub)}{\mp\bm\mt}
\frac{\mp\eta(\bm\mt,\ub)}{\mp\bm\mt\TT} \, \xi(d\ub) \,.
\end{equation}
The equality $\Cb(w,\xi,\mtbb)=\Mb^{-1}(\xi,\mtbb)$ is obtained
for $w(u)=c\, \ms^{-2}(u)$, with $c$ a positive constant, and this
choice of $w(\cdot)$ is thus optimal (in terms of asymptotic
variance) among all WLS estimators. This result can be compared to
that obtained for linear regression in Section \ref{S:weighing}
where $\ms^2\Mb_N^{-1}$ was the {\em exact expression} for the
variance of $\mthb^N$ for $N$ finite. In nonlinear regression the
expression $\Cb(w,\xi,\mtbb)/N$ for the variance of $\mthb^N$ is
only valid asymptotically, see (\ref{ANWLS}); moreover, it depends
on the unknown true value $\mtbb$ of the parameters. These results
can easily be extended to situations where also the variance of
the errors depends on the parameters $\bm\mt$ of the response
$\eta$, that is,
$\Ex_{\ub_k}(\mve_k^2)=\ms^2(\ub_k)=\beta\ml(\mtbb,\ub_k)$, see
e.g.\ \cite{PPmoda72004}, \cite{PP-eusipco-04}.

\subsection{Maximum-likelihood estimation} \label{S:ML}

Denote $\varphi_{\ub_k}(\cdot)$ the probability density function
(pdf) of the error $\mve_k$ in (\ref{regression}). Due to the
independence of errors, we obtain for the vector $\yb$ of
observation the pdf $\pi(\yb|\bm\mt)=\prod_{k=1}^N
\pi[y(\ub_i)|\bm\mt] = \prod_{k=1}^N \varphi_{\ub_k}
[y(\ub_k)-\eta(\bm\mt,\ub_k)]$ and the Maximum-Likelihood (ML)
estimator $\mthb_{ML}^N$ minimizes
$-\log\pi(\yb|\bm\mt)=\sum_{k=1}^N
-\log\varphi_{\ub_k}[y(\ub_k)-\eta(\bm\mt,\ub_k)]$. Different pdf
$\varphi$ yield different estimators (LS for Gaussian errors,
$L_1$ estimation for errors with a Laplace distribution, etc.).
Under standard regularity assumptions on $\varphi_\ub(\cdot)$ and
for designs satisfying conditions (i) or (ii) of Section
\ref{S:WLS}, $\mthb_{ML}^N\raas\mtbb$ and
\begin{equation}\label{ANML}
    \sqrt{N} (\mthb^N_{ML}-\mtbb) \rad z
\sim\SN(0,\Mb_F^{-1}(\xi,\mtbb))\,, \ N\ra\infty \,,
\end{equation}
with $\Mb_F(\bm\mt,\xi)$ the Fisher information matrix (average
per sample) given by
\begin{eqnarray*}
  \Mb_F(\bm\mt,\xi) &=& \Ex_\mt\left\{ \frac{1}{N}
\frac{\mp\log\pi(\yb|\bm\mt)}{\mp\bm\mt}
\frac{\mp\log\pi(\yb|\bm\mt)}{\mp\bm\mt\TT} \right\} \\
  &=& -\Ex_\mt\left\{ \frac{1}{N}
\frac{\mp^2\log\pi(\yb|\bm\mt)}{\mp\bm\mt\mp\bm\mt\TT} \right\}
\,.
\end{eqnarray*}
In the particular case of the regression model considered here we
obtain
\begin{equation}\label{MF1}
        \Mb_F(\xi,\bm\mt)= \int_\SU \SI(\ub) \,
\frac{\mp\eta(\bm\mt,\ub)}{\mp\bm\mt}
\frac{\mp\eta(\bm\mt,\ub)}{\mp\bm\mt\TT} \, \xi(d\ub)
\end{equation}
with $\SI(\ub)=\int [\varphi'_\ub(z)]^2/\varphi_\ub(z) \, dz$ the
Fisher information for location of the pdf $\varphi_\ub$. From the
Cram\'er-Rao inequality, $\Mb_F^{-1}(\xi,\bm\mt)$ forms a
lower-bound on the covariance matrix of any unbiased estimator
$\mthb^N$ of $\bm\mt$, i.e.,
$\Ex_\mt\{(\mthb^N-\bm\mt)(\mthb^N-\bm\mt)\TT\} -
\Mb_F^{-1}(\xi,\bm\mt)/N$ is non-negative definite for any
estimator $\mthb^N$ such that $\Ex_\mt\{\mthb^N\}=\bm\mt$. When
the errors $\mve_k$ are normal $\SN(0,\ms^2(\ub_k))$,
$\SI(\ub)=\ms^{-2}(\ub)$ and ML estimation coincides with WLS with
optimal weights (and $\Mb_F(\xi,\bm\mt)$ coincides with
(\ref{MinfoWLS})). When they are i.i.d., that is
$\varphi_\ub=\varphi$ for any $\ub$, $\SI(\ub)=\SI$ constant, and
\begin{equation}\label{MinfoF}
    \Mb_F(\xi,\bm\mt)= \SI \int_\SU \frac{\mp\eta(\bm\mt,\ub)}{\mp\bm\mt}
\frac{\mp\eta(\bm\mt,\ub)}{\mp\bm\mt\TT} \, \xi(d\ub) \,.
\end{equation}

\subsubsection{Estimating functions}\label{S:Estimatinge}

Estimating functions form a very generally applicable set of tools
for parameter estimation in stochastic models. As the example
below will illustrate, they can yield very simple estimators for
dynamical systems. One can refer to \cite{Heyde97} for a general
exposition of the methodology, see also the discussion paper
\cite{LiangZ95} that comprises a short historical perspective.
Instrumental variables methods (see, e.g., \cite{SoderstromS81},
\cite{SoderstromS83} and Chapter 8 of \cite{SoderstromS89}) used
in dynamical systems as an alternative to LS estimation when the
regressors and errors are correlated (so that the LS estimator is
biased) can be considered as methods for constructing unbiased
estimating functions. Their implementation often involves the
construction of regressors obtained through simulations with
previous values of parameter estimates, but simpler constructions
are possible.

Consider a discrete-time system with scalar state and input,
respectively $x_i$ and $u_i$, defined by the recurrence equation
\begin{equation}\label{discrete}
    x_{i+1}=x_i + T[u_i+\mtb(x_i+1)]\,, \ i=0,1,2\ldots
\end{equation}
with known sampling period $T$ and initial state $x_0$. The
observations are given by $y_i=x_i+\mve_i$ for $i\geq 1$, where
$(\mve_i)$ denotes a sequence of i.i.d.\ errors normal
$\SN(0,\ms^2)$. The unknown parameter $\mtb$ can be estimated by
LS (which corresponds to ML estimation since the errors are
normal), but recursive LS cannot be used since $x_i$ depends
nonlinearly in $\mtb$.  However, simpler estimators can be used if
one is prepared to loose some precision for the estimation. For
instance, substitute $y_i$ for the state $x_i$ in (\ref{discrete})
and form the equation in $\mt$
\begin{equation}\label{gi}
    g_{i+1}(\mt)= y_{i+1} - y_i - T[u_i+\mt(y_i+1)] = 0 \,;
\end{equation}
$k$ successive observations then give $G_k(\mt) = \frac1k
\sum_{i=1}^{k} g_i(\mt)$ $= 0$. Since $G_k(\mt)$ is linear in the
$y_i$'s, $\Ex_\mt\{G_k(\mt)\}=0$ for any $\mt$, and $G_k(\mt)$ is
called an {\em unbiased estimating
function}\label{PFN:NL-EE}\footnote{\label{FN:NL-EE}Nonlinearity
in the observations is allowed, provided that the bias is suitably
corrected; for instance the function $g'_{i+1}(\mt)=(1+y_i)
g_{i+1}(\mt) + \ms^2(1+T\mt)$ with $g_{i+1}(\mt)$ given by
(\ref{gi}) satisfies $\Ex_\mt\{g'_{i+1}(\mt)\}=0$ for any $\mt$
when the errors $\mve_i$ are i.i.d.\ with zero mean and variance
$\ms^2$, and $(1/k)\sum_{i=1}^{k} g'_i(\mt)$ is an unbiased
estimating function for $\mt$.}, see, e.g., \cite{LiangZ95}. Since
$G_k(\mt)$ is linear in $\mt$, the solution $\tilde\mt^k$ of
$G_k(\mt)=0$ is simply given by
\begin{equation}\label{tildemt}
    \tilde\mt^{k} = \frac{(y_{k}-y_0)/(kT) - (\sum_{i=0}^{k-1}
u_i)/k}{1+(\sum_{i=0}^{k-1} y_i)/k}
\end{equation}
(provided that the denominator is different from zero) and forms
an estimator for $\mtb$. Notice that the true value $\mtb$
satisfies a similar equation with the $y_i$'s replaced by the
noise-free values $x_i$. Estimation by $\tilde\mt^{k}$ is less
precise than LS estimation, see Figure \ref{F:KKK-thetas} in
Section \ref{S:further-issues}, but requires much less
computations. Would other parameters be present in the model,
other estimating functions would be required. For instance, a
function of the type $G_{k,\ma}(\bm\mt)=\sum_{i=1}^{k} i^\ma
g_i(\bm\mt)$ would put more stress on the transient (respectively
long-term) behavior of the system when $\ma<0$ (respectively
$\ma>0$). Also, the multiplication of $g_i(\bm\mt)$ by a known
function of $u_i$ gives a new estimating function. When
information on the noise statistics is available, it is desirable
for the (asymptotic) precision of the estimation to choose $G_k$
as (proportional to) an approximation of the score function $\mp
\log \pi(\yb|\bm\mt)/\mp\bm\mt$ with $\pi(\yb|\bm\mt)$ the pdf of
the observations $y_1,\ldots,y_k$, see, e.g., \cite{CoxH74} p.\
274 and \cite{LiangZ95}.

There seems to be a revival of interest for estimating functions,
partly due to the elegant algebraic framework recently developed
for time-continuous linear systems (differential equations); see
\cite{FliessS-R2003} where estimating functions are constructed
²through Laplace transforms. However, in this algebraic setting
only multiplications by $s$ or $s^{-1}$ and differentiation with
respect to $s$ are considered (with $s$ the Laplace variable),
which seems unnecessarily restrictive. Consider for instance the
time-continuous version of (\ref{discrete}),
\begin{equation}\label{diff-eq}
 \dot x = u +  \mtb (x+1) \,, \ x(0)=x_0 \,,
\end{equation}
where $\dot x$ denotes differentiation with respect to time. Its
Laplace transform is $s X(s) = U(s)+ \mtb X(s) + s^{-1} \mtb +
x_0$, which can be first multiplied by $s$, then differentiated
two times with respect to $s$ and the result multiplied by
$s^{-2}$ to avoid derivation with respect to time. This gives a
estimating function comprising double integrations with respect to
time. Multiple integrations may be avoided by noticing that the
multiplication of the initial differential equation by any
function of time preserves the linearity of the estimating
function with respect to both $\mt$ and the state (provided that
the integrals involved are well defined). For instance, when $\dot
u$ is a known function of time, the multiplication of
(\ref{diff-eq}) by the input $u$ followed by integration with
respect to time gives the estimating function
$[x(t)u(t)-x_0u_0]/t- (1/t) \int_0^t [x(\tau)\dot
u(\tau)+u^2(\tau)]\, d\tau = (\mt/t) \, \int_0^t
u(\tau)[1+x(\tau)] \, d\tau$, which is linear in $x$. Infinitely
many unbiased estimating functions can thus be easily constructed
in this way. (Note that, due to linearity, the introduction of
process noise in (\ref{diff-eq}) as $\dot x(t) = u(t) +
\mtb[x(t)+1] + dB_t(\omega)$, with $B_t(\omega)$ a Brownian
motion, leaves the estimating function above unbiased.)

The analysis of the asymptotic behavior of the estimator $\mttb^k$
associated with an estimating function is straightforward when the
function is unbiased and linear in $\bm\mt$. The expression of the
asymptotic variance of the estimator can be used to select
suitable experiments in terms of the precision of the estimation,
as it is the case for LS or ML estimation. However, in general the
asymptotic variance of the estimator takes a more complicated form
than $\Mb^{-1}(\xi,\bm\mt)$ or $\Mb_F^{-1}(\xi,\bm\mt)$, see
(\ref{MinfoWLS}, \ref{MinfoF}), so that DOE for such estimators
does not seem to have been considered so far. The recent revival
of interest for this method might provide some motivation for such
developments (see also Section \ref{S:further-issues}).


\subsection{DOE}

To obtain a precise estimation of $\bm\mt$ one should first use a
good estimator (WLS with weights proportional to $\ms^{-2}$, or
ML) and second select a good design\footnote{We shall thus follow
the standard approach, in which the estimator is chosen first, and
an optimal design is then constructed for that given estimator
(even though it may be optimal for different estimators); this can
be justified under rather general conditions, see
\cite{Nather75}.} $\xi^*$. In the next section we shall consider
classical DOE for parameter estimation, which is based on the
information matrix (\ref{MinfoF})\footnote{Note that defining
$\tilde\eta(\bm\mt,\ub)=\ms^{-1}(\ub)\eta(\bm\mt,\ub)$ and
$\tilde\eta(\bm\mt,\ub)=\sqrt{\SI(\ub)}\eta(\bm\mt,\ub)$ one can
respectively write the matrices (\ref{MinfoWLS}) and (\ref{MF1})
in the same form as (\ref{MinfoF}). Also notice that classical DOE
uses the covariance matrix with the simplest expression: DOE for
WLS estimation is more complicated for non-optimal weights than
for the optimal ones, compare $\Cb(w,\xi,\bm\mt)$ to
$\Mb^{-1}(\xi,\bm\mt)$ in Section \ref{S:WLS}. Similarly, the
asymptotic covariance matrix for a general $M$-estimator (see,
e.g., \cite{Huber81}) is more complicated than for ML.}. Hence, we
shall choose $\xi^*$ that optimizes $\Phi[\Mb_F(\xi,\bm\mt)]$, for
some criterion function $\Phi(\cdot)$. For models nonlinear in
$\mt$, this raises two difficulties: (i) the criterion function,
and thus $\xi^*$, depends on a guessed value $\bm\mt$ for $\mtbb$.
This is called {\em local} DOE (the design $\xi^*$ is optimal
locally, when $\mtbb$ is close to $\bm\mt$), some alternatives to
local optimal design will be presented in Section
\ref{S:design-nonlinear}; (ii) the method relies on the {\em
asymptotic properties} of the estimator. More accurate
approximations of the precision of the estimation exist, see e.g.\
\cite{Pazman93}, but are complicated and seldom used for DOE, see
\cite{PazmanPa92}, \cite{PPa94} (see also the recent work
\cite{CampiW2005} concerning the finite sample size properties of
estimators, which raises challenging DOE issues). They will not be
considered here. For dynamical systems with correlated
observations or containing an autoregressive part, classical DOE
also relies on the information matrix, which has then a more
complicated expression, see Section \ref{S:control-for-design}.
Also, the calculation of the asymptotic covariance of some
estimators requires specific developments that are not presented
here, see e.g. \cite{GoodwinP77}, \cite{Ljung87}, \cite{Caines88}
for recursive estimation methods. For Bayesian estimation, a
standard approach for DOE consists in replacing
$\Mb_F(\xi,\bm\mt)$ by $\Mb_F(\xi,\bm\mt)+\Omega^{-1}/N$, with
$\Omega$ the prior covariance matrix for $\bm\mt$, see e.g.\
\cite{Pilz83}, \cite{ChalonerV95}. Note finally the central role
of the design concerning the asymptotic properties of estimators.
In particular, the conditions (i) and (ii) of Section \ref{S:WLS}
on the design imply some stationarity of the ``inputs'' $\ub_k$
and guarantee the {\em persistence of excitation}, which can be
expressed as a condition on the minimum eigenvalue of the
information matrix: $\liminf_{N\ra\infty} \ml_{\min}
[\Mb_F(\xi_N,\bm\mt)] >0$, with $\xi_N$ the empirical measure of
$\ub_1,\ldots,\ub_N$ (that is, $\liminf_{N\ra\infty}
\ml_{\min}(\Mb_N) /N >0$ for the linear regression model of
Section \ref{S:weighing}, see (\ref{Mlinear})).

\section{DOE for parameter estimation} \label{S:design-for-estimation}

\subsection{Design criteria}\label{S:criteria}

We consider criteria for designing optimal experiments (for
parameter estimation) that are scalar functions of the (Fisher)
information matrix (average, per sample)
(\ref{MinfoF})\footnote{Notice that the analytic form of the
sensitivities $\mp \eta(\bm\mt,\ub)/\mp\bm\mt$ of the model
response is not required: for a model given by differential
equations, like in Section \ref{S:pharmaco}, or by difference
equations, the sensitivities can be obtained by simulation,
together with the model response itself; see, e.g., Chapter 4 of
\cite{WPl97}.}. For $N$ observations at the design points
$\ub_i\in\SU$, $i=1,\ldots,N$, we shall denote
$U_1^N=(\ub_1,\ldots,\ub_N)$, which is called a {\em finite} (or
{\em discrete}) design of size $N$, or {\em $N$-point design}. The
associated information matrix is then
\begin{equation}\label{MinfoFN}
    \Mb_F(U_1^N, \bm\mt)= \frac{\SI}{N} \sum_{i=1}^N
\frac{\mp\eta(\bm\mt,\ub_i)}{\mp\bm\mt}
\frac{\mp\eta(\bm\mt,\ub_i)}{\mp\bm\mt\TT} \,.
\end{equation}
The admissible design set $\SU$ is sometimes a finite set,
$\SU=\{\ub^1,\ldots,\ub^K\}$, $K<\infty$. We shall more generally
assume that $\SU$ is a compact subset of $\mathbb{R}^d$. For a
linear regression model with i.i.d.\ errors $\SN(0,\ms^2)$, the
ellipsoid $\SR(\mthb_{LS}^N,\ma)=\{\bm\mt \ / \
(\bm\mt-\mthb_{LS}^N)\TT\Mb_F(U_1^N)(\bm\mt-\mthb_{LS}^N)\leq
\SX_\ma^2(p)/N \}$, where $\SX_\ma^2(p)$ has the probability $\ma$
to be exceeded by a random variable chi-square distributed with
$p$ degrees of freedom, satisfies
$\Pr\{\mtbb\in\SR(\mthb_{LS}^N,\ma)\}=\ma$, and this is
asymptotically true in nonlinear situations\footnote{Such
confidence regions for $\bm\mt$ can be transformed into
simultaneous confidence regions for functions of $\bm\mt$, see in
particular \cite{Scheffe61}, \cite{BomboisAG2005}.}.

Most of classical design criteria are related to characteristics
of (asymptotic) confidence ellipsoids. Minimizing
$\Phi(\Mb)=\tr[\Mb^{-1}]$ corresponds to minimizing the sum of the
squared lengthes of the axes of (asymptotic) confidence ellipsoids
for $\bm\mt$ and is called $A$-optimal design (minimizing
$\Phi(\Mb)=\tr[\Qb\TT\Qb\Mb^{-1}]$ with $\Qb$ some weighting
matrix is called $L$-optimal design, see \cite{Chernoff53} for an
early reference). Minimizing the longest axis of (asymptotic)
confidence ellipsoids for $\bm\mt$ is equivalent to maximizing the
minimum eigenvalue of $\Mb$ and is called $E$-optimal design.
$D$-optimal design maximizes $\det(\Mb)$, or equivalently
minimizes the volume of (asymptotic) confidence ellipsoids for
$\bm\mt$ (their volume being proportional to $1/\sqrt{\det\Mb}$).
This approach is very much used, in particular due to the
invariance of a $D$-optimal experiment by re-parametrization of
the model (since $\det\Mb(\xi,\bm\mt') = \det\Mb(\xi,\bm\mt)
[\det(\mp\bm\mt'/\mp\bm\mt\TT)]^{-2}$). Most often $D$-optimal
experiments consist of replications of a small number of different
experimental conditions. This has been illustrated by the example
of Section \ref{S:pharmaco} for which $p=4$ and four sampling
times were duplicated in the $D$-optimal design $\tb^*$.

\subsection{Algorithms for discrete design}\label{S:algorithms-exact}

Consider the regression model (\ref{regression}) with i.i.d.\
errors and $N$ observations at $U_1^N=(\ub_1,\ldots,\ub_N)$ where
the {\em support points} $\ub_i$ belong to
$\SU\subset\mathbb{R}^d$. The Fisher information matrix
$\Mb_F(U_1^N,\bm\mt)$ is then given by (\ref{MinfoFN}). The
(local) design problem consists in optimizing
$\Psi_\mt(U_1^N)=\Phi[\Mb_F(U_1^N,\bm\mt)]$ for a given $\bm\mt$,
with respect to $U_1^N \in \mathbb{R}^{N\times d}$. If the problem
dimension $N\times d$ is not too large, standard optimization
algorithms can be used (note, however, that constraints may exist
in the definition of the admissible set $\SU$ and that local
optima exist is general). When $N \times d$ is large, specific
algorithms are recommended. They are usually of the exchange type,
see \cite{Fedorov72}, \cite{Mitchell74}. Since several local
optima exist in general, these methods provide locally optimal
solutions only.

\subsection{Approximate design theory} \label{S:approximate-designs}

\subsubsection{Design measures}\label{S:construction}

Suppose that replications of observations exist, so that several
$\ub_i$'s coincide in (\ref{MinfoFN}). Let $m<N$ denote the number
of different $\ub_i$'s, so that
$$
\Mb_F(U_1^N, \bm\mt)= \SI \sum_{i=1}^m \frac{r_i}{N}
\frac{\mp\eta(\bm\mt,\ub_i)}{\mp\bm\mt}
\frac{\mp\eta(\bm\mt,\ub_i)}{\mp\bm\mt\TT}
$$
with $r_i/N$ the proportion of observations collected at $\ub_i$,
which can be considered as the {\em percentage of experimental
effort} at $\ub_i$, or the weight of the support point $\ub_i$.
Denote $\ml(\ub_i)$ this weight. The design $U_1^N$ is then
characterized by the support points $\ub_1,\ldots,\ub_m$ and their
associated weights $\ml(\ub_1), \ldots, \ml(\ub_m)$ satisfying
$\sum_{i=1}^m \ml(\ub_i)=1$, that is, a normalized discrete
distribution on the $\ub_i$'s, with the constraints
$\ml(\ub_i)=r_i/N$, $i=1,\ldots,m$. Releasing these constraints,
one defines an {\em approximate design} as a discrete probability
measure with support points $\ub_i$ and weights $\ml_i$ (with
$\sum_{i=1}^m \ml_i=1$). Releasing now the discreteness
constraint, a {\em design measure} is simply defined as any
probability measure $\xi$ on $\SU$, see \cite{KieferW59}, and
$\Mb_F(\xi,\bm\mt)$ takes the form (\ref{MinfoF}). Now,
$\Mb_F(\xi,\bm\mt)$ belongs to the convex hull of the set $\SM_1$
of rank-one matrices of the form $\Mb(\delta_\ub,\bm\mt) = \SI \,
[\mp\eta(\bm\mt,\ub)/\mp\bm\mt] \,
[\mp\eta(\bm\mt,\ub)/\mp\bm\mt\TT]$. It is a $p\times p$ symmetric
matrix, and thus belongs to a $p(p+1)/2$-dimensional space.
Therefore, from Caratheodory's Theorem, it can be written as the
linear combination of $p(p+1)/2+1$ elements of $\SM_1$ at most;
that is
\begin{equation}\label{MinfoFdiscrete}
    \Mb_F(\xi,\bm\mt) = \SI \sum_{i=1}^m \ml_i
\frac{\mp\eta(\bm\mt,\ub_i)}{\mp\bm\mt}
\frac{\mp\eta(\bm\mt,\ub_i)}{\mp\bm\mt\TT} \,,
\end{equation}
with $m \leq p(p+1)/2 + 1$. The information matrix associated with
any design measure $\xi$ can thus always be considered as obtained
from a discrete probability measure with $p(p+1)/2+1$ support
points at most. This is true in particular for the optimal
design\footnote{In general the situation is even more favorable.
For instance, if $\xi_D$ is $D$-optimal (it maximizes $\det
\Mb_F(\xi,\bm\mt)$), then $\Mb_F(\xi_D,\bm\mt)$ is on the boundary
of the convex closure of $\SM_1$ and $p(p+1)/2$ support points are
enough.}. Given such a discrete design measure $\xi$ with $m$
support points, a discrete design $U_1^N$ with repetitions can be
obtained by choosing the numbers of repetitions $r_i$ such that
$r_i/N$ is an approximation\footnote{This is at the origin of the
name {\em approximate design theory}. However, a design $\xi$
(even with a density) can sometimes be implemented {\em without
any approximation}: this is the case in Section
\ref{S:frequency-domain} where $\xi$ corresponds to the power
spectral density of the input signal.} of $\ml_i$, the weight of
$\ub_i$ for $\xi$, see, e.g., \cite{PukelsheimR92}.

The property that the matrices in the sum (\ref{MinfoFdiscrete})
have rank one is not fundamental here and is only due to the fact
that we considered single-output models (i.e., scalar
observations). In the multiple-output case with independent
errors, say with $\yb(\ub)$ of dimension $q$ corrupted by errors
having the $q\times q$ covariance matrix $\bm\Sigma(\ub)$, the
model response is a $q$-dimensional vector $\bm\eta(\bm\mt,\ub)$
and the information matrix for WLS estimation with weights
$\bm\Sigma^{-1}(\ub)$ is $\Mb(\xi,\bm\mt)= \int_\SU
[\mp\bm\eta\TT(\bm\mt,\ub)/\mp\bm\mt]\, \bm\Sigma^{-1}(\ub) \,
[\mp\bm\eta(\bm\mt,\ub)/\mp\bm\mt\TT] \, \xi(d\ub)$, to be
compared with (\ref{MinfoWLS}) obtained in the single-output case,
see, e.g., \cite{Fedorov72}, Section 1.7 and Chapter 5.
Caratheodory's Theorem still applies and, with the same notations
as above, we can write
$$
\Mb(\xi,\bm\mt) = \sum_{i=1}^m \ml_i \,
\frac{\mp\bm\eta\TT(\bm\mt,\ub_i)}{\mp\bm\mt}
\bm\Sigma^{-1}(\ub_i)
\frac{\mp\bm\eta(\bm\mt,\ub_i)}{\mp\bm\mt\TT} \,,
$$
with again $m \leq p(p+1)/2 + 1$. All the results concerning DOE
for scalar observations thus easily generalize to the
multiple-output situation.

\subsubsection{Properties}\label{S:properties- approximateD}

Only the main properties are indicated, one may refer to the books
\cite{Fedorov72}, \cite{Silvey80}, \cite{Pazman86},
\cite{AtkinsonD92}, \cite{Pukelsheim93}, \cite{FedorovH97} for
more detailed developments. Suppose that the design criterion
$\Phi[\Mb]$ to be minimized (respectively maximized) is strictly
convex (respectively concave). For instance for $D$-optimality,
maximizing $\det[\Mb]$ is equivalent to maximizing $\log\det[\Mb]$
and, for any positive-definite matrices $\Mb_1$, $\Mb_2$ such that
$\Mb_1\neq \Mb_2$, $\forall \alpha$, $0<\alpha<1$, $\log
\det[(1-\alpha)\Mb_1+\alpha \Mb_2] > (1-\alpha)\log \det [\Mb_1] +
\alpha \log \det [\Mb_2]$, so that $\Phi[\cdot]= \log \det[\cdot]$
is a strictly concave function. Since $\Mb_F(\xi,\bm\mt)$ belongs
to a convex set, the optimal matrix $\Mb_F^*=\Mb_F(\xi^*,\bm\mt)$
for $\Phi$ is unique (which usually does not imply that the
optimal design $\xi^*$ is unique; however, the set of optimal
design measures is convex). The uniqueness of the optimum and
differentiability of the criterion directly yield {\em a necessary
and sufficient condition for optimality}, and in the case of
$D$-optimality we obtain the following, known as {\em
Kiefer-Wolfowitz Equivalence Theorem} \cite{KieferW60} (other
equivalence theorems are easily obtained for other design criteria
having suitable regularity and the appropriate convexity or
concavity property).

\begin{theo}\label{Th:KWET} The following statements are equivalent:\\
(1) $\xi_D$ is $D$-optimal for $\bm\mt$, \\
(2) $\max_{\ub\in\SU} d_\mt(\ub,\xi_D) = p$, \\
(3) $\xi_D$ minimizes $\max_{\ub\in\SU} d_\mt(\ub,\xi_D)$,\\
where $d_\mt(\ub,\xi)$ is defined by
\begin{equation}\label{d}
 d_\mt(\ub,\xi) = \SI \,
\frac{\mp\eta(\bm\mt,\ub)}{\mp\bm\mt\TT} \, \Mb_F^{-1}(\xi,\bm\mt)
\, \frac{\mp\eta(\bm\mt,\ub)}{\mp\bm\mt} \,.
\end{equation}
Moreover, for any support
point $\ub_i$ of $\xi_D$, $d_\mt(\ub_i,\xi_D)=p$.
\end{theo}

Note that condition (2) is easily checked when $u$ is scalar by
plotting $d_\mt(u,\xi)$ as a function of $u$.

Theorem \ref{Th:KWET} relates optimality in the parameter space to
optimality in the space of observations, in the following sense.
Let $\mthb_{ML}^N$ be obtained for a design $\xi$, the variance of
the prediction $\eta(\mthb_{ML}^N,\ub)$ of the response at $\ub$
is then such that $N \var[\eta(\mthb_{ML}^N,\ub)]$ tends to
\begin{equation}\label{var}
    \frac{\mp\eta(\bm\mt,\ub)}{\mp\bm\mt\TT}_{|\mtbb} \,
    \Mb_F^{-1}(\xi,\mtbb)\,
\frac{\mp\eta(\bm\mt,\ub)}{\mp\bm\mt}_{|\mtbb} =
\frac{d_{\mtb}(\ub,\xi)}{\SI}
\end{equation}
when $N\ra\infty$, see (\ref{ANML}). Therefore, a $D$-optimal
experiment also minimizes the maximum of the (asymptotic) variance
of the prediction over the experimental domain $\SU$. This is
called $G$-optimality, and Theorem \ref{Th:KWET} thus expresses
the {\em equivalence between $D$ and $G$-optimality}. (It is also
related to maximum entropy sampling considered in Section
\ref{S:Model-based}, see \cite{Wynn2004}.)

Suppose that the observations are collected sequentially and that
the choice of the design points can be made accordingly ({\em
sequential design}). After the collection of
$y(\ub_1),\ldots,y(\ub_N)$, which gives the parameter estimates
$\mthb^N$ and the prediction $\eta(\mthb^N,\ub)$, in order to
improve the precision of the prediction the next observation
should intuitively be placed where $\var[\eta(\mthb^N,\ub)]$ is
large, that is, where $d_{\mth^N}(\ub,\xi_N)$ is large, with
$\xi_N$ the empirical measure for the first $N$ design points.
This receives a theoretical justification in the algorithms
presented below.

\subsubsection{Algorithms}\label{S:algorithms-for-approximateD}

The presentation is for $D$-optimality, but most algorithms easily
generalize to other criteria. Let $\xi^k$ denote the design
measure at iteration $k$ of the algorithm. The steepest-ascent
direction at $\xi^k$ corresponds to the delta measure that puts
mass 1 at $\ub^*_{k+1}=\arg\max_{\ub\in\SU} d_\mt(\ub,\xi^k)$.
Hence, at iteration $k$, algorithms of the steepest-ascent type
add the support point $\ub^*_{k}$ to $\xi^k$ as follows:

{\bf Fedorov--Wynn Algorithm:} \\
$\bullet$ Step 1~: Choose $\xi^1$ not degenerate ($\det
\Mb_F(\xi^1,\bm\mt)\neq 0$), and $\me$ such that $0<\me<<1$, set $k=1$. \\
$\bullet$ Step 2~: Compute $\ub_{k+1}^* = \arg\max_{\ub\in\SU}
d_\mt(\ub,\xi^k)$. If $d_\mt(\ub_{k+1}^*,\xi^k) < p+\me$, stop:
$\xi^k$ is almost $D$-optimal. \\
$\bullet$ Step 3~: Set $\xi^{k+1}=(1-\alpha_k)\xi^k+\alpha_k
\delta_{\ub_{k+1}^*}$, $k\ra k+1$, return to Step 2.

Fedorov's algorithm corresponds to choosing the step-length
$\alpha_k^*$ that maximizes $\log\det \Mb_F(\xi^{k+1},\bm\mt)$,
which gives $\alpha_k^*=
[d_\mt(\ub_{k+1}^*,\xi^k)-p]/\{p[d_\mt(\ub_{k+1}^*,\xi^k)-1]\}$
(note that $0<\alpha_k^*<1/p$) and ensures monotonic convergence
towards a $D$-optimal measure $\xi_D$, see \cite{Fedorov72}.

Wynn's algorithm corresponds to a sequence satisfying
$0<\alpha_k<1$, $\lim_{k\ra\infty} \alpha_k=0$ and
$\sum_{i=1}^\infty \alpha_k = \infty$, see \cite{Wynn70} (the
convergence is then not monotonic). One may notice that in
sequential design where the design points enter
$\Mb_F(U_1^N,\bm\mt)$ given by (\ref{MinfoFN}) one at a time, one
has
\begin{eqnarray*}
\Mb_F(U_1^{k+1},\bm\mt) &=& \frac{k}{k+1} \, \Mb_F(U_1^k,\bm\mt) \\
&& + \frac{1}{k+1} \, \SI\,
\frac{\mp\eta(\bm\mt,\ub_{k+1})}{\mp\bm\mt}
\frac{\mp\eta(\bm\mt,\ub_{k+1})}{\mp\bm\mt\TT}
\end{eqnarray*}
and, when $\ub_{k+1}=\arg\max_{\ub\in\SU} d_\mt(\ub,\xi^k)$, this
corresponds to Wynn's algorithm with $\alpha_k=1/(k+1)$.

Contrary to the exchange algorithms of Section
\ref{S:algorithms-exact}, these steepest-ascent methods guarantee
convergence to the optimum. However, in practice they are rather
slow (in particular due to the fact that a support point present
at iteration $k$ is never totally removed in subsequent iterations
---~since $\alpha_k<1$ for any $k$) and faster methods, still of
the steepest-ascent type, have been proposed, see e.g.
\cite{Bohning89}, \cite{MolchanovZ01}, \cite{MolchanovZ02} and
\cite{FedorovH97} p.\ 49. An acceleration of the algorithms can
also be obtained by using a submodularity property of the design
criterion, see \cite{RobertazziS89}, or by removing design points
that cannot support a $D$-optimal design measure, see
\cite{HPa06_SPL}.

When the set $\SU$ is finite (which can be obtained by a suitable
discretization), say with cardinality $K$, the optimal design
problem in the approximate design framework corresponds to the
minimization of a convex function of $K$ positive weights $\ml_i$
with sum equal one, and any convex optimization algorithm can be
used. The recent progress in interior point methods, see for
instance the survey \cite{ForsgrenGW2002} and the books
\cite{NesterovN94}, \cite{denHertog94}, \cite{Wright97},
\cite{Ye97}, provide alternatives to the usual sequential
quadratic programming algorithm. In control theory these methods
have lead to the development of tools based on linear matrix
inequalities, see, e.g., \cite{BoydEGFB94}, which in turn have
been suggested for $D$-optimal design, see \cite{VandenbergheBW98}
and Chapter 7 of \cite{BoydV2004}. Alternatively, a simple
updating rule can sometimes be used for the optimization of a
design criterion over a finite set $\SU=\{\ub^1,\ldots,\ub^K\}$.
For instance, convergence to a $D$-optimal measure is guaranteed
when the weight $\ml_i^k$ of $\ub^i$ at iteration $k$ is updated
as
\begin{equation}\label{Ben}
\ml_i^{k+1} = \ml_i^k \frac{ d_\mt(\ub^{i},\xi^k) }{p} \,,
\end{equation}
where $\xi^k$ is the measure defined by the support points $\ub^i$
and their associated weights $\ml_i^k$, and $d_\mt(\ub,\xi)$ is
given by (\ref{d}), see \cite{Titterington76}, \cite{SilveyTT78},
\cite{Torsney83} and Chapter 5 of \cite{Pazman86}. (Note that
$\sum_{i=1}^K \ml_i^{k+1}=1$ and that $\ml_i^{k+1}>0$ when
$\ml_i^k>0$.) The extension to the case where information matrices
associated with single points have ranks larger than one (see
Section \ref{S:construction}) is considered in \cite{Ucinski2005}.

Finally, it is worthwhile noticing that $D$-optimal design is
connected with a minimum-ellipsoid problem. Indeed, using
Lagrangian theory one can easily show that the construction of
$\xi_D$ that maximizes the determinant of $\Mb_F(\xi,\bm\mt)$
given by (\ref{MinfoF}) with respect to the probability measure
$\xi$ on $\SU$ is equivalent to the construction of the
minimum-volume ellipsoid, centered at the origin, that contains
the set $\SS_\mt=\{\mp\eta(\bm\mt,\ub)/\mp\bm\mt: \, \ub\in\SU\}
\subset\mathbb{R}^p$, see \cite{Sibson72}. The construction of the
minimum-volume ellipsoid centered at 0 containing a given set
$\SU\subset\mathbb{R}^p$ thus corresponds to a $D$-optimal design
problem on $\SU$ for the linear regression model
$\eta(\bm\mt,\ub)=\ub\TT\bm\mt$. In the case where the center of
the ellipsoid is free, one can show equivalence with a $D$-optimal
design in a $(p+1)$-dimensional space where the regression model
is $\eta(\bm\mt,\ub)=(1 \ \ \ub\TT)\bm\mt$,
$\bm\mt\in\mathbb{R}^{p+1}$, see \cite{SilvermanT80},
\cite{Titterington75}. Algorithms with iterations of the type
(\ref{Ben}) are then strongly connected with steepest-descent type
algorithms when minimizing a quadratic function, see
\cite{PWZa01-AAM}, \cite{PWZa05_MP} and Chapter 7 of
\cite{PWZl2000}. In system identification, minimum-volume
ellipsoids find applications in parameter bounding (or parameter
estimation with bounded errors), see, e.g., \cite{PWa94_2} and
\cite{RaynaudPW01a} for an application to robust control.

\subsubsection{Active learning with parametric
models}\label{S:active-parametric}

When learning with a parametric model, the prediction
$\yh_\SD(\ub)$ at $\ub$ is $\eta(\mthb^N,\ub)$ with $\mthb^N$
estimated from the data
$\SD=\{[\ub_1,y(\ub_1)],\ldots,[\ub_N,y(\ub_N)]\}$. As Theorem
\ref{Th:KWET} shows, the precision of the prediction is directly
related to the precision of the estimation of the model parameters
$\bm\mt$: a $D$-optimal design minimizes the maximum (asymptotic)
variance\footnote{We could also speak of MSE since in parametric
models the estimators are usually unbiased for models linear in
$\bm\mt$, and for nonlinear models (under the condition of
persistence of excitation) the squared bias decreases as $1/N^2$
whereas the variance decreases as $1/N$, see \cite{Box71}.} of
$\yh_\SD(\ub)$ for $\ub\in\SU$. Similar properties hold for other
measures of the precision of the prediction. Consider for instance
the integrated (asymptotic) variance of the prediction with
respect to some given probability measure $\pi$ (that may express
the importance given to different values of $\ub$ in $\SU$). It is
given by $\Psi_{\mt,\Hb}(\xi)=\tr \left\{ \Hb \Mb^{-1}(\xi,\bm\mt)
\right\}$, where $\Hb = \Hb(\bm\mt) = \int_{\SU}
[\mp\eta(\bm\mt,\ub)/\mp\bm\mt]\,
[\mp\eta(\bm\mt,\ub)/\mp\bm\mt\TT] \, \pi(d\ub)$, see (\ref{var}),
and its minimization corresponds to a $L$-optimal design problem,
see Section \ref{S:criteria}. The following parametric learning
problem is addressed in \cite{Kanamori2002}: the measure $\pi$ is
unknown, $n$ samples $\ub_i$ from $\pi$ are used, together with
the associated observations, to estimate $\bm\mt$ and $\Hb$,
respectively by $\mthb^n$ and $\hat\Hb^n(\mthb^n)$, $N-n$ samples
are then chosen optimally for $\Psi_{\mth^n,\hat\Hb^n}(\xi)$. It
is shown that the optimal balance between the two sample sizes
corresponds to $n$ being proportional to $\sqrt{N}$. When the
samples $\ub_i$ are cheap and only the observations $y(\ub_i)$ are
expensive, one may decide on-line to collect an observation or not
for updating the estimate $\mthb^n$ and the information matrix
$\Mb_n$. A sequential selection rule is proposed in \cite{Pa05},
which is asymptotically optimal when a given proportion $n=
\lfloor \ma N \rfloor$ of samples, $\ma\in(0,1)$, can be accepted
in a sequence of length $N$, $N\ra\infty$.

There exists a fundamental difference between learning with
parametric and nonparametric models. For parametric models, the
MSE of the prediction globally decreases as $1/N$, and precise
predictions are obtained for optimal designs which, from
Caratheodory's Theorem (see Section \ref{S:construction}) are
concentrated on a finite number of sites. These are the points
$\ub_i$ that carry the maximum information about $\bm\mt$ useful
for prediction, in terms of the selected design criterion. On the
opposite, precise predictions for nonparametric models are
obtained when the observation sites are spread over $\SU$, see
Section \ref{S:Model-based}. Note, however, that {\em parametric
methods rely on the extremely strong assumption that the data are
generated by a model with known structure}. Since optimal designs
will tend to {\em repeat observations at the same sites} (whatever
the method used for their construction), {\em modelling errors
will not be detected}. This makes optimal design theory of very
delicate use when the model is of the behavioral type, e.g.\ a
neural network as in \cite{Cohn94}, \cite{CohnGJ96}. A recent
approach \cite{GazutMI2006} based on bagging (Bootstrap
Aggregating, see \cite{Breiman96}) seems to open promising
perspectives.

\subsubsection{Dependence in $\bm\mt$ in nonlinear
situations}\label{S:design-nonlinear}

We already stressed the point that in nonlinear situations the
Fisher information matrix depends on $\bm\mt$, so that an optimal
design for estimation depends on the unknown value of the
parameters to be estimated. So far, only {\em local optimal
design} has been considered, where the experiment is designed for
a nominal value $\bm\mt$. Several methods can be used to reduce
the effect of the dependence in the assumed $\bm\mt$. A first
simple approach is to use a finite set
$\Theta=\{\bm\mt^{(1)},\ldots,\bm\mt^{(m)}\}$ of nominal values
and to design $m$ locally optimal experiments $\xi^*_{\mt^{(i)}}$
for the $\bm\mt^{(i)}$'s in $\Theta$. This permits to appreciate
the strength of the dependence of the optimal experiment in
$\bm\mt$, and several $\xi^*_{\mt^{(i)}}$'s can eventually be
combined to form a single experiment. More sophisticated
approaches rely on average or minimax optimality.

In {\em average-optimal design}, the criterion
$\Psi_\mt(\xi)=\Phi[\Mb_F(\xi,\bm\mt)]$ is replaced by its
expectation $\Ex_\pi\{ \Psi_\mt(\xi) \} = \int
\Phi[\Mb_F(\xi,\bm\mt)]\, \pi(d\bm\mt)$ for some suitably chosen
prior $\pi$, see, e.g., \cite{Fedorov80}, \cite{ChalonerL89},
\cite{ChalonerV95}. (Note that when the Fisher information matrix
$\Mb_F(\xi,\bm\mt)$ is used, it means that the prior is not used
for estimation and the method is not really Bayesian.) In {\em
minimax-optimal design}, $\Psi_\mt(\xi)$ (to be minimized) is
replaced by its worst possible value $\max_{\bm\mt\in\Theta}
\Phi[\Mb_F(\xi,\mt)]$ when $\bm\mt$ belongs to a given feasible
set $\Theta$, see, e.g., \cite{Fedorov80}. Compared to local
design, these approaches do not create any special difficulty
(other than heavier computations) for {\em discrete design}, see
Section \ref{S:algorithms-exact}: no special property of the
design criterion is used, but the algorithms only yield local
optima. Of course, for computational reasons the situation is
simpler when $\pi$ is a discrete measure and $\Theta$ is a finite
set\footnote{When $\Theta$ is a compact set of $\mathbb{R}^p$, a
relaxation algorithm is suggested in \cite{PWa88} for
minimax-optimal design; stochastic approximation can be used for
average-optimal design, see \cite{PWa85}.}. Concerning {\em
approximate design theory} (Section \ref{S:approximate-designs}),
the convexity (or concavity) of $\Phi$ is preserved, Equivalence
Theorems can still be obtained (Section \ref{S:properties-
approximateD}) and globally convergent algorithms can be
constructed (Section \ref{S:algorithms-for-approximateD}), see,
e.g., \cite{FedorovH97}. A noticeable difference with local
design, however, concerns the number of support points of the
optimum design which is no longer bounded by $p(p+1)/2+1$ (see,
e.g., Appendix A in \cite{RojasWGF07}). Also, algorithms for
minimax-optimal design are more complicated than for local optimal
design, in particular since the steepest-ascent direction does not
necessarily correspond to a one-point delta measure.

A third possible approach to circumvent the dependence in $\bm\mt$
consists in designing the experiment sequentially (see the
examples in Sections \ref{S:discrimination} and
\ref{S:optimization}), which is particularly well suited for
nonparametric models, both in terms of prediction and estimation
of the model, see Section \ref{S:Model-based}. Sequential DOE for
regression models is considered into more details in Section
\ref{S:sequential-design}.

\section{Control in DOE: optimal inputs for parameter estimation
in dynamical models} \label{S:control-for-design}

In this section, the choice of the input is (part of) the design,
$U_1^N$ or $\xi$ depending whether discrete or approximate design
is used. One can refer in particular to the book \cite{Zarrop79}
and Chapter 6 of \cite{GoodwinP77} for detailed developments. The
presentation is for single-input single-output systems, but the
results can be extended to multi-input multi-output systems. The
attention is on the construction of the Fisher information matrix,
the inverse of which corresponds to the asymptotic covariance of
the ML estimator, see Section \ref{S:information-matrix}. For
control-oriented applications it is important to relate the
experimental design criterion to the ultimate control objective,
see, e.g., \cite{ForssellL2000}, \cite{Gevers2005}. This is
considered in Section \ref{S:frequency-domain}.

\subsection{Input design in the time domain}\label{S:time-domain}

Consider a Box and Jenkins model, with observations
$$
y_k=F(\mtbb,z)u_k+G(\mtbb,z)\mve_k
$$
where the errors $\mve_k$ are i.i.d.\ $\SN(0,\ms^2)$, and
$F(\bm\mt,z)$ and $G(\bm\mt,z)$ are rational fractions in $z^{-1}$
with $G$ stable with a stable inverse. Suppose that $\ms^2$ is
unknown. An extended vector of parameters $\bm\beta=(\bm\mt\TT \ \
\ms^2)\TT$ must then be estimated, and one can assume that
$G(\bm\mt,\infty)=1$ without any loss of generality. For suitable
input sequences (such that the experiment is informative enough,
see \cite{Ljung87}, p.\ 361), $N \Var(\bm{\hat\beta}_{ML}^N) \ra
\Mb_F^{-1}(\xi,\bm{\bar\beta})$, $N\ra\infty$, with
$\bm{\bar\beta}$ the unknown true value of $\bm\beta$ and
$$
\Mb_F(\xi,\bm\beta) = \Ex_{\beta}\left\{\frac{1}{N} \frac{\mp \log
\pi(\yb|\bm\beta)}{\mp\bm\beta} \frac{\mp \log
\pi(\yb|\bm\beta)}{\mp\bm\beta\TT} \right\} \,.
$$
Using the independence and normality of the errors and the fact
that $\ms^2$ does not depend on $\bm\mt$, we obtain
$$
\Mb_F(\xi,\bm\beta) = \left( \begin{array}{cc}  \Mb_F(\xi,\bm\mt) & \0b \\
\0b\TT & \frac{1}{2\ms^4} \end{array}\right)
$$
$$
\mbox{with } \Mb_F(\xi,\bm\mt) = \Ex_\mt\left\{\frac{1}{N\ms^2}
\sum_{k=1}^N \frac{\mp e_k(\bm\mt)}{\mp\bm\mt}\frac{\mp
e_k(\bm\mt)}{\mp\bm\mt\TT} \right\}
$$
and $e_k(\bm\mt)$ the prediction error
$e_k(\bm\mt)=G^{-1}(\bm\mt,z)[y_k-F(\bm\mt,z)u_k]$. The fact that
$\ms^2$ is unknown has therefore no influence on the (asymptotic)
precision of the estimation of $\bm\mt$. Assuming that the
identification is performed in open loop (that is, there is no
feedback)\footnote{One may refer, e.g., to Chapter 6 of
\cite{GoodwinP77}, \cite{GeversL85}, \cite{GeversL86},
\cite{HjalmarssonGD-B96}, \cite{ForssellL99}, \cite{ForssellL2000}
\cite{JanssonH-IFAC2005} for results concerning closed-loop
experiments.} and that $F$ and $G$ have no common parameters (that
is, $\bm\mt$ can be partitioned into $\bm\mt=(\bm\mt_F\TT \ \
\bm\mt_G\TT)\TT$, with $p_F$ components in $\bm\mt_F$ and $p_G$ in
$\bm\mt_G$), we then obtain
$$
    \Mb_F(\xi,\bm\mt) = \left( \begin{array}{cc}  \Mb_F^F(\xi,\bm\mt) & \Ob
\\ \Ob &  \Mb_F^G(\xi,\bm\mt) \end{array}\right)
$$
with
\begin{eqnarray}
\Mb_F^F(\xi,\bm\mt)  &=& \frac{1}{N\ms^2} \sum_{k=1}^N
\left[ G^{-1}(\bm\mt,z) \frac{\mp F(\bm\mt,z)}{\mp\bm\mt_F} u_k \right] \nonumber \\
&&  \hspace{1.5cm} \times \left[ G^{-1}(\bm\mt,z) \frac{\mp
F(\bm\mt,z)}{\mp\bm\mt_F\TT} u_k \right]  \label{MFF}
\end{eqnarray}
and $\Mb_F^G(\xi,\bm\mt)$ not depending on $\{u_k\}$, see, e.g.,
\cite{GoodwinP77}, p.\ 131. The asymptotic covariance matrix
$\Mb_F^{-1}(\xi,\bm\mt)$ is thus partitioned into two blocks, and
the input sequence $(u_k)$ has no effect on the precision of the
estimation of the parameters $\bm\mt_G$ in $G$. A $D$-optimal
input sequence maximizes $\det \Mb_F^F(\xi,\bm\mt) = \det\left[
1/(N\ms^2) \sum_{k=1}^N \vb_k\vb\TT_k\right]$ where $\vb_k$ is a
vector of (linearly) filtered inputs,
\begin{equation}\label{vk}
    \vb_k=G^{-1}(\bm\mt,z) \frac{\mp F(\bm\mt,z)}{\mp\bm\mt_F} \, u_k \,,
\end{equation}
usually with power or amplitude constraints on $u_k$. This
corresponds to an optimal control problem in the time domain and
standard techniques from control theory can be used for its
solution.

\subsection{Input design in the frequency
domain}\label{S:frequency-domain}

We consider the same framework as in previous section, with the
information matrix of interest $\Mb_F^F(\xi,\bm\mt)$ given by
(\ref{MFF}). Suppose that the system output is uniformly sampled
at period $T$ and denote
$\underline{\Mb}_F^F(\xi,\bm\mt)=\lim_{N\ra\infty}
\Mb_F^F(\xi,\bm\mt)/T$ the average Fisher information matrix per
time unit. It can be written as $\underline{\Mb}_F^F(\xi,\bm\mt) =
1/(2\pi\ms^2) \int_{-\pi}^{\pi} \SP_v(\omega) \, d\omega$ with
$\SP_v(\omega)$ the power spectral density of $\vb_k$ given by
(\ref{vk}), or $\underline{\Mb}_F^F(\xi,\bm\mt) = 1/\pi
\int_{0}^{\pi} \tilde{\underline{\Mb}}_F^F(\omega,\bm\mt) \,
\SP_u(\omega) \,d\omega$ with $\SP_u(\omega)$ the power spectral
density of $u$ and
\begin{eqnarray*}
\tilde{\underline{\Mb}}_F^F(\omega,\bm\mt) &=& \frac{1}{\ms^2}
\SR_e \left\{ \frac{\mp
F(\bm\mt,e^{j\omega})}{\mp\bm\mt_F} G^{-1}(\bm\mt,e^{j\omega})  \right.\\
&& \hspace{1cm} \times \left. G^{-1}(\bm\mt,e^{-j\omega})
\frac{\mp F(\bm\mt,e^{-j\omega})}{\mp\bm\mt_F\TT} \right\} \,.
\end{eqnarray*}
The framework is thus the the same as for approximate design
theory of Section \ref{S:approximate-designs}: the experimental
domain $\SU$ becomes the frequency domain $\mathbb{R}^+$ and to
the design measure $\xi$ corresponds the power spectral density
$\SP_u$. An optimal input with discrete spectrum always exists; it
has a finite number of support points\footnote{One can show that
the upper bound on their number can be reduced from $p_F(p_F+1)/2$
to $p_F$, the number of parameters in $F$, see \cite{GoodwinP77},
p.\ 138.} (frequencies) and associated weights (input power). The
optimal input can thus be searched in the class of signals
consisting of finite combinations of sinusoidal components, and
the algorithms for its construction are identical to those of
Section \ref{S:algorithms-for-approximateD}. Notice, however, that
no approximation is now involved in the implementation of the
``approximate'' design. Once an optimal spectrum has been
specified, the construction of signal with this spectrum can obey
practical considerations, for instance on the amplitude of the
signal, see \cite{BarkerG99}. Alternatively, the input spectrum
can be decomposed on a suitable basis of rational transfer
functions and the optimization of $\SP_u$ performed with respect
to the linear coefficients of the decomposition, see
\cite{JanssonH2004}, \cite{JanssonH2005}. Notice that the problem
can also be taken the other way round: one may wish to minimize
the input power subject to a constraint on the precision of the
estimation, expressed through $\Mb_F^{-1}(\xi,\bm\mt)$, see e.g.,
\cite{BomboisSGHVdH2004}, \cite{BomboisSGVdHH2005}.

The design criteria presented in Section \ref{S:criteria} are
related to the definition of confidence regions, or uncertainty
sets, for the model parameters. When the intended application of
the identification is the control of a dynamical system, it seems
advisable {\em to relate the DOE to control-oriented uncertainty
sets}, see in particular \cite{Gevers2005} for an inspired
exposition. First note that according to the expression
(\ref{var}) the variance of the transfer function
$F(\bm\mt,e^{j\omega})$ at the frequency $\omega$ is approximately
$V_F(\omega) = (1/N) [\mp
F(\bm\mt,e^{j\omega})/\mp\bm\mt_F\TT]\,[\Mb_F^F(\xi,\bm\mt)]^{-1}\,[\mp
F(\bm\mt,e^{j\omega})/\mp\bm\mt_F]$. Several $H_\infty$-related
design criteria can then be derived. For instance, a
robust-control constraint of the form $\|W(e^{j\omega}) \Delta
F(\bm\mt,e^{j\omega})/F(\bm\mt,e^{j\omega})\|_\infty <1$, with
$\Delta F(\bm\mt,e^{j\omega})/F(\bm\mt,e^{j\omega})$ the relative
error on $F(\bm\mt,e^{j\omega})$ due to the estimation of $\bm\mt$
and $W(e^{j\omega})$ a weighting function, leads to
$\zb\TT(\bm\mt,e^{j\omega})[\Mb_F^F(\xi,\bm\mt)]^{-1}\zb(\bm\mt,e^{j\omega})
<1$ $\forall \omega$, with $\zb(\bm\mt,e^{j\omega}) =
(1/\sqrt{N})|W(e^{j\omega})|[\mp
F(\bm\mt,e^{j\omega})/\mp\bm\mt_F]$. This type of constraint can
be expressed as a linear matrix inequality in
$\Mb_F^F(\xi,\bm\mt)$, and, using the KYP lemma, the problem can
be reformulated as having a finite number of constraints, see
\cite{JanssonH2005}. Notice that minimizing $\max_\omega
\zb\TT(\bm\mt,e^{j\omega})[\Mb_F^F(\xi,\bm\mt)]^{-1}\zb(\bm\mt,e^{j\omega})$
can be compared to $E$-optimum design, see Section
\ref{S:criteria}, which minimizes $\max_{\{\zb: \zb\TT\zb=1\}}
\zb\TT[\Mb_F^F(\xi,\bm\mt)]^{-1}\zb$. When $|W(e^{j\omega})|=1$
(uniform weighting) and $G(\bm\mt,e^{j\omega})=1$ (white noise),
it corresponds to $G$-optimal design, and thus to $D$-optimal
design, see Section \ref{S:properties- approximateD}. It is also
strongly related to the minimax-optimal design of Section
\ref{S:design-nonlinear}, (where the worst-case is now considered
with respect to $\omega$), see \cite{FedorovH97} and \cite{PWa88}
for algorithms. Alternatively, the asymptotic confidence regions
for $\bm\mt$ can be transformed into uncertainty sets
$\SS_F(\bm\mt,\xi)$ for the transfer function
$F(\bm\mt,e^{j\omega})$. The worst-case $\nu$-gap over this set
can then be computed, with the property that the smaller this
number, the larger the set of controllers that stabilize all
transfer functions in $\SS_F(\bm\mt,\xi)$ \cite{GeversBCSA2003a},
\cite{GeversBCSA2003b} (see also \cite{RaynaudPW01a} for related
results). Designing experiments that minimize the worst-case
$\nu$-gap is considered in \cite{HildebrandG2003} where the
problem is shown to be amenable to convex optimization.

{\em The dependence of the design criteria in the unknown
parameters of the model is a major issue for optimal input
design}, as it is more generally the case for models with a
nonlinear parametrization (it explains why input spectra with few
sinusoidal components are often considered as unpleasant). The
methods suggested in Section \ref{S:design-nonlinear} to face this
difficulty can be applied here too. In particular, input spectra
having a small number of components can be avoided by designing
optimal inputs for different nominal values for $\bm\mt$ and
combining the optimal spectra that are obtained, or by using
average or minimax-optimal design \cite{RojasWGF07}. One can also
design the experiment sequentially (see Section
\ref{S:sequential-design}); in general, each design step involves
many observations and a few steps only are required to achieve
suitable performance, see, e.g. \cite{BarenthinJH2005}.

When on-line adaptation is possible, adjusting the controller
while data are collected and the uncertainty on the model
decreases can be expected to achieve better performance than
non-adaptive robust control. Ideally, one would wish to have
uncertainty sets shrinking towards a single point representing the
true model (or the model closest to the true system for the model
class considered), so that a robust controller adapted to smaller
and smaller uncertainty sets becomes less and less conservative.
While the determination of such robust-and-adaptive controllers is
still an open issue, a first step in the construction is to
investigate the properties of the parameter estimates in adaptive
procedures.

\section{DOE in adaptive control} \label{S:adaptive-control}

The results of Sections \ref{S:design-for-estimation} and
\ref{S:control-for-design} rely on the asymptotic properties of
the estimator: the asymptotic variance of $\mthb_{ML}^N$ was
supposed to be given by $\Mb_F^{-1}/N$, which is true when the
design (input) sequence satisfies some ``stationarity'' condition
(the assumption of {\em random design} was used in Section
\ref{S:design-for-estimation} and a condition of {\em persistence
of excitation} in Section \ref{S:control-for-design}). However,
this condition may fail to hold: a typical example is adaptive
control, where the input has another objective than estimation.
The issues that it raises are investigated hereafter. We first
present a series of simple examples that illustrate the variety of
the difficulties.

\subsection{Examples of difficulties}\label{S:exLai&Wei}

It is rather well-known that the usual asymptotic normality of the
LS estimator may fail to hold for designs such that $\Mb_F(U_1^N)$
is nonsingular for any $N$ but converges to a singular matrix,
that is, such that $\ml_{\min}[\Mb_F(U_1^N)] \ra 0$ as
$N\ra\infty$, see \cite{PPa05_SPL}. We shall not develop this
point but rather focuss on the difficulties raised by the
sequential construction of the design.

Consider the following well-known example (see, e.g.,
\cite{LaiR81}, \cite{LaiW82}) of a linear regression model with
observations $y_k=\mtb_1 + \mtb_2 u_k + \mve_k$ where the errors
$\mve_k$ are i.i.d.\ with zero mean and variance 1. The input
(design points $u_k$) satisfies $u_1=0$ and
$u_{n+1}=(1/n)\sum_{i=1}^n u_i+ (c/n) \sum_{i=1}^n \mve_i$. Then,
one can prove that $\{\mthb_{LS}^N\}_1 \raas \mtb_1+
\sum_{i=1}^\infty \mve_i/i$ and $\{\mthb_{LS}^N\}_2 \raas \mtb_2
-1/c$, $N\ra\infty$. That is, $\{\mthb_{LS}^N\}_1$ converges to a
random variable and $\{\mthb_{LS}^N\}_2$ to a non-random constant
different from $\mtb_2$. The non-consistency of the LS estimator
is due to the dependence of $u_{n+1}$ on previous $\mve_i$'s, that
is, to the presence of feedback control (in terms of DOE, the
design is sequential). Although $\Mb_N=\sum_{i=1}^N (1 \ \
u_i)\TT(1 \ \ u_i)$ is such that $\ml_{\min}(\Mb_N)\ra \infty$, it
does not grows fast enough (in particular, the information matrix
$\Mb(U_1^N)=\Mb_N/N$ tends to become singular). Although this
example might seem quite artificial, one must notice that adaptive
control as used e.g.\ in self-tuning strategies, may raise similar
difficulties.

\subsubsection{ARX model and self-tuning regulator}\label{S:STR}
Consider a model with observations satisfying $y_k=a_1
y_{k-1}+\cdots+a_{n_a} y_{k-n_a} +b_1
u_{k-1}+\cdots+b_{n_b}u_{k-n_b}+\mve_k$, which we can write
$y_k=\rb_k\TT \mtbb + \mve_k$, with
$\mtbb=(b_1,b_2,\ldots,a_1,a_2,\ldots)\TT$ and
$\rb_k=(u_{k-1},\ldots,u_{k-n_b},$ $y_{k-1},\ldots,y_{k-n_a})\TT$.
The objective of minimum-variance control is to minimize
$R_N=\sum_{k=1}^N (y_k-\mve_k)^2$. The input sequence $(u_k)$ is
then said to be {\em globally convergent} if $R_N/N \raas 0$ as
$N\ra\infty$, see \cite{LaiW86}, \cite{LaiW87}, \cite{Guo94}. If
$\mtbb$ is known (with $b_1\neq 0$) the optimal controller
corresponds to $u^*_{k}=-(a_1 y_{k}+\cdots+a_{n_a} y_{k+1-n_a}+b_2
u_{k-1}+\cdots+b_{n_b}u_{k+1-n_b})/b_1$. But then
$\rb_k\TT\mtbb=0$ for all $k$, the matrix $\Mb_N=\sum_{k=1}^N
\rb_k \rb_k\TT$ is singular (since $\mtbb\TT \Mb_N \mtbb=0$) and
$\mtbb$ is not estimable. If certainty equivalence is forced by
using at step $k$ the optimal control calculated for
$\mthb_{LS}^k$, then additional perturbations must be introduced
to guarantee that $\ml_{\min}(\Mb_N)$ tends to infinity fast
enough, see, e.g., \cite{AstromW73}. Using a persistently exciting
input $u_k$, possibly with optimal features via the approach of
Section \ref{S:control-for-design}, permits to avoid this
difficulty but is in conflict with the global convergence property
\cite{LaiW86}, in particular since $\|\mtbb\|^2 \ml_{\min}(\Mb_N)
< R_N$, see \cite{Guo94}.

\subsubsection{Self-tuning optimizer}\label{S:STO}

Consider a linear regression model with observations
$y_k=\rb\TT(\ub_k) \mtbb+\mve_k$. The objective is to maximize a
function $f(\ub,\mtbb)$ with respect to $\ub$. If $\mtbb$ were
known, the value $\ub^*=\ub^*(\mtbb)=\arg\max_\ub f(\ub,\mtbb)$
could be used (for instance, $u^*=-\mtb_1/(2\mtb_2)$ when
$f(u,\bm\mt)=\mt_0+\mt_1 u + \mt_2 u^2$). Since $\mtbb$ is
unknown, it must be estimated from the observations $y_k$,
$k=1,2,\ldots$ Again, the matrix $\Mb_N=\sum_{k=1}^N \rb(\ub_k)
\rb\TT(\ub_k)$ is singular when the control is fixed, that is when
$\ub_k=\ub^*(\bm\mt)$ (constant) for all $k$, and $\bm\mt$ is then
not estimable. Suppose that forced certainty equivalence is used
with LS estimation, that is $\ub_{k+1}=\ub^*(\mthb^k_{LS})$.
Perturbations should then be introduced to ensure consistency
(e.g. randomly, see \cite{BozinZ91} for the quadratic case
$f(u,\bm\mt)=\mt_0+\mt_1 u + \mt_2 u^2$). The persistency of
excitation is here in conflict with the performance objective
$(1/n)\sum_{i=1}^n f(\ub_i,\mtbb) \raas f(\ub^*,\mtbb)$,
$n\ra\infty$. Self-tuning regulation of dynamical systems is
considered in \cite{KrsticW2000} and \cite{Krstic2000} for
time-continuous systems and in \cite{ChoiKAL2002} for
discrete-time systems. With a periodic disturbance of magnitude
$\ma$ playing the role of a persistently exciting input signal,
the output exponentially converges to a neighborhood $\SO(\ma^2)$
of the extremum.

\subsection{Nonlinear feedback control is not the
answer}\label{S:NFC}

Nonlinear-Feedback Control (NFC) offers a set of techniques for
stabilizing systems with unknown parameters, see in particular the
book \cite{KrsticKK95}. The stability of the closed-loop is proved
using Lyapunov techniques and, although not explicitly expressed
in the construction of the feedback control, an estimator of the
model parameters is obtained, which differs from standard
estimation methods. At first sight one might think that NFC brings
a suitable answer to adaptive control issues. However, stability
is not consistency and it is the aim of this section to show that
a direct application of NFC is bound to fail in the presence of
random disturbances. Combining NFC with more traditional
estimation methods and suitably exciting perturbations then forms
interesting perspectives, see Section \ref{S:further-issues}.

\begin{figure}[h]
\begin{center}
\includegraphics[scale=0.35]{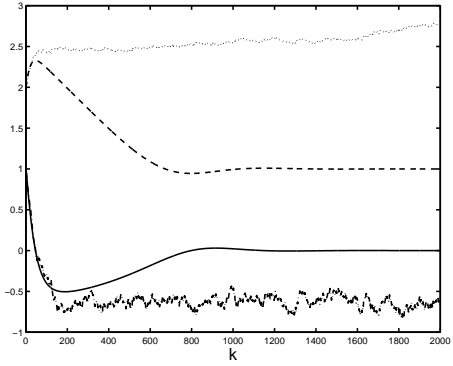}
\caption{Evolution of $x_k$ (solid line) and $\mth_k$ (dashed
line) as functions of $k$ for the system (\ref{discrete}) with NFC
(\ref{NFC-nonoise}) ($\mtb=1$, $\mth_0=2$, $x_0=1$, $a=1$,
sampling period $T=0.01$ s). The curves in dash-dotted line and
dotted line respectively show $x_k$ and $\mth_k$ when $y_k$ is
substituted for $x_k$ in (\ref{NFC-nonoise}).}
\label{F:KKK-nonoise}
\end{center}
\end{figure}

The presentation is made through (a slight modification of) one of
the simplest examples in \cite{KrsticKK95}. Consider the dynamical
system (\ref{diff-eq}), with known initial state $x_0$ and unknown
parameter $\mtb\in\mathbb{R}$. The problem is to construct a
control $u=u(t)$ that drives $x$ to zero. (Notice that if $\mtb$
were known, $u=-(a+\mtb)x-\mtb$ with $a>0$ would solve the problem
since substitution in (\ref{diff-eq}) gives the stable system
$\dot x=-ax$.) The following method is suggested in
\cite{KrsticKK95}: (i) construct an auxiliary controller that
obeys $\dot{\mth} = x(x+1)$, (ii) consider $\mth$ as an estimator
of $\mtb$ and use FCE control with $\mth$, that is, $u=-(a+
\mth)x- \mth$, $a>0$. The stability of this NFC can be checked
through the behavior of the Lyapunov function $V(x,\mth)=x^2/2+
(\mt-\mth)^2/2$. It satisfies $\dot V(x,\mth)=-ax^2$, which
implies that $x$ tends to zero, as required. Then, $\mth+u$ tends
to zero (from the expression of $u$), and $\mtb+u$ also tends to
zero (from Lasalle principle). Therefore, the estimation error
$\mth-\mtb$ tends to zero\footnote{In the example on p.\ 3--4 of
\cite{KrsticKK95}, $\dot x = u + \mtb x$, $\dot{\mth} = x^2$ and
$u=-(a+ \mth)x$, $a>0$, so that $x$ tends to zero but not
necessarily the estimation error $\mth-\mtb$.}. In the simulations
that follow we simply use a discretized Euler approximation of the
differential equation (\ref{diff-eq}) and of the associated
continuous-time controller, although it should be emphasized that
some care is needed in general when implementing a digital
controller on a continuous-time model, see, e.g.,
\cite{NesicT2001}. The discretization of (\ref{diff-eq}) gives the
recurrence equation (\ref{discrete}), where $x_k=x(kT)$ and
$u_k=u(kT)$. We take $\mtb=1$ and $x_0=1$, the sampling period $T$
is taken equal to $0.01$~s. (Notice that the open-loop system is
unstable.) The NFC is discretized as
\begin{equation}
\label{NFC-nonoise}
\begin{array}{l}
\mth_{k+1} = \mth_k + T x_k(x_k+1)\,, \\ u_k = -(a+ \mth_k)x_k-
\mth_k \,,
\end{array}
\end{equation}
where $\mth_k=\mth(kT)$. We take $\mth_0=2$ and $a=1$ s$^{-1}$.
(Although the book \cite{KrsticKK95} only concerns the
stabilization of continuous-time systems, one can easily check
that the fixed point $x_k=0$, $\mth_k=\mtb$ of the controlled
discretized model above is Lyapunov-asymptotically stable.)
Simulation results are presented in Figure \ref{F:KKK-nonoise}.
The initial decrease of the state variable (solid line) is in
agreement with the time-constant $a^{-1}=1$~s and, for $t>8$~s,
the parameter estimates (dashed line) and state become very close
to the targets, respectively $\mtb=1$ and zero.

Suppose now that the state is observed through $y_0=x_0$ and
$y_k=x_k+\mve_k$ for $k\geq 1$, where $(\mve_k)$ denotes a
sequence of i.i.d.\ errors normal $\SN(0,\ms^2)$ (setting
$\sigma^2=S^2T$ one may suppose for instance that $\mve_k$ is $S$
times the integral of a realization of the standard Brownian
motion between 0 and $T$). We take $\ms =x_0/2=0.5$, a rather
extreme situation, to emphasize the influence of measurement
errors. The evolutions of $x_k$ (dash-dotted line) and $\mth_k$
(dotted line) when $y_k$ is substituted for $x_k$ in
(\ref{NFC-nonoise}) are  presented on Figure \ref{F:KKK-nonoise}:
the sequence of parameter estimates does not converge, the state
fluctuates and is clearly not driven to zero.

\subsection{Some consistency results}\label{S:some-results}

The difficulties encountered in Sections \ref{S:STR}, \ref{S:STO}
and \ref{S:NFC} are general in regulation-type problems: in order
to satisfy the control objective, the input should asymptotically
vanish, which does not bring enough excitation for guaranteeing
the consistent estimation of the model parameters. The control
objective is thus in conflict with parameter estimation, and
perturbations must be introduced. It is then of importance to know
the minimal amount of perturbations required to ensure consistency
of the estimator on which the control is based. Some results are
presented below for the case of linear regression.

\subsubsection{LS estimation}\label{S:results-LS}

Consider a linear regression model with observations
$y_k=\rb_k\TT\mtbb  + \mve_k$, and denote by $\SF_k$ the
$\ms$-algebra generated by the errors $\mve_1,\ldots,\mve_k$. They
are supposed to form a martingale difference sequence ($\mve_k$ is
$\SF_{k-1}$ measurable and $\Ex\{\mve_k|\SF_{k-1}\}=0$) and to be
such that $\sup_k \Ex\{\mve_k^2|\SF_{k-1}\}<\infty$ (with i.i.d.\
errors with zero mean and finite variance as a special case). Let
$\Mb_N=\sum_{k=1}^N \rb_k \rb_k\TT$, then $\Mb_N^{-1} \ra 0$ for
$N\ra\infty$ is \\
$\bullet$  sufficient for $\mthb_{LS}^N \raas \mtbb$ when the
regressors $\rb_k$ are non-random constants, see \cite{LaiRW78},
\cite{LaiRW79};\\
$\bullet$ necessary and sufficient if, moreover, the errors are
$\mve_k$ i.i.d.,\\
$\bullet$ but $\Mb_N^{-1} \raas 0$ is not sufficient for
$\mthb_{LS}^N \raas \mtbb$ if $\rb_k$ is $\SF_{k-1}$ measurable
(see the first example of Section \ref{S:exLai&Wei}).

In the latter situation, a sufficient condition for $\mthb_{LS}^N
\raas \mtbb$ when $N\ra\infty$ is that $\ml_{\min}(\Mb_N) \raas
\infty $ and $ [\log \ml_{\max}(\Mb_N)]^{1+\delta} =
o[\ml_{\min}(\Mb_N)]$ a.s.\ for some $\delta>0$, see
\cite{LaiW82}. In some sense, this is the best possible condition:
it is only marginally violated in the first example of Section
\ref{S:exLai&Wei}, where $[\log
\ml_{\max}(\Mb_N)]/\ml_{\min}(\Mb_N)$ tends a.s.\ to a random
constant. Note that this condition is much weaker than the
persistence of excitation which requires that $\Mb_N$ grows at the
same speed as $N$.

\subsubsection{Bayesian imbedding}\label{S:results-bayesian}

An even weaker condition is obtained for Bayesian estimation. Let
$\pi$ be a prior probability measure for $\bm\mt$ and denote by
$P$ the probability measure induced by the errors $\mve_k$,
$k=1,\ldots,\infty$. Denote $\SF'_k$ the $\ms$-algebra generated
by the observations $y_1,\ldots,y_k$ and suppose that $\rb_k$ is
$\SF'_{k-1}$-measurable. Suppose that the parameters are estimated
by the posterior mean $\mthb_{B}^N = \Ex\{\bm\mt|\SF'_N\}$ and
denote by $\Cb_N=\Var(\bm\mt|\SF'_N)$ the posterior covariance
matrix. Then, from martingale theory, $\mthb_{B}^N$ and $\Cb_N$
both converge $(\pi\times P)$-a.s.\ when $N\ra\infty$, see
\cite{Sternby77}, and all what is required for the $(\pi\times
P)$-a.s.\ consistency of the estimator is $\Cb_N\ra 0$ $(\pi\times
P)$-a.s. Now, for a linear regression model with i.i.d.\ normal
errors $\mve_k$ and a normal prior for $\bm\mt$, Bayesian
estimation coincides with LS estimation (when the prior for
$\bm\mt$ is suitably chosen), $\Cb_N$ is proportional to
$\Mb_N^{-1}$ and therefore, $\Mb_N^{-1} \ra 0$ $(\pi\times
P)$-a.s.\ is sufficient for $\mthb_{LS}^N=\mthb_B^N \ra \mtbb$
$(\pi\times P)$-a.s. The required condition is thus as weak as
when the regressors $\rb_k$ are non-random constants! Note,
however, that the convergence is almost sure with respect to
$\mtbb$ having the prior $\pi$; that is, singular values of
$\mtbb$ may exist for which consistency does not hold\footnote{In
the first example of Section \ref{S:exLai&Wei} $\rb_{k}$ is not
$\SF'_{k-1}$-measurable since $u_k$ is not obtained from previous
observations. Modify the control into $u_{n+1}= \ma_1 +
(\ma_2/n)\sum_{i=1}^n u_i + (c/n) \sum_{i=1}^n y_i$, which is
$\SF'_{n}$-measurable. Then, $\mthb_{LS}^N$ is not consistent when
$\mtbb$ takes the particular value $\mtb_1=-\ma_1/c$,
$\mtb_2=(1-\ma_2)/c$, so that the control coincides with previous
one, $u_{n+1}=(1/n)\sum_{i=1}^n u_i+ (c/n) \sum_{i=1}^n \mve_i$.}.

This very powerful technique which analyses the properties of LS
estimation via a Bayesian approach is called {\em Bayesian
imbedding}, see \cite{Sternby77}, \cite{Kumar90}. Although in its
original formulation it requires the measurement errors to be
normal, the normality assumption is relaxed in \cite{Hu96} to the
condition that the density $\varphi$ of the errors is log-concave
($d^2 \log \varphi(t)/dt^2<0$) and strictly positive with respect
to the Lebesgue measure $\mu$, the prior measure $\pi$ being
absolutely continuous with respect to $\mu$. More generally, the
consistency of Bayesian estimators can be checked through the
behavior of posterior covariance matrices, see \cite{Hu97}.
Bayesian imbedding allows for easier proofs of consistency of the
estimator, and permits to relax the conditions on the
perturbations required to obtain consistency. This is illustrated
below by revisiting the examples of Sections \ref{S:STR} and
\ref{S:STO}.

Consider again the self-tuning regulator of Section \ref{S:STR}.
When LS estimation is used with forced certainty equivalence
control, it is required to perturb the system to obtain a globally
convergent input. It can be shown \cite{Lai86} that the control
objective $R_n$ grows at least as $\log n$, and randomly perturbed
input sequences achieving this performance are proposed in
\cite{LaiW87}. Using Bayesian imbedding, global convergence can be
obtained without the introduction of perturbations, see
\cite{Kumar90}.

For the self-tuning optimizer of Section \ref{S:STO}, {\AA}str\"om
and Wittenmark \cite{AstromW89} have suggested a control of the
type $\ub_{k+1}=\arg\max_{\ub} f(\ub,\mthb^k) +  \ma_k \
d(\ub,\xi_k)/k$, where
$d(\ub,\xi)=\rb\TT(\ub)\Mb^{-1}(\xi)\rb(\ub)$ is the function
(\ref{d}) used in $D$-optimal design, $\xi_k$ is the empirical
measure of the inputs $\ub_1,\ldots,\ub_k$ and $(\ma_k)$ is a
sequence of positive scalars. Note that $d(\ub,\xi_k)/k=
\rb\TT(\ub)\Mb_k^{-1}\rb(\ub)$ with $\Mb_k=\sum_{i=1}^k \rb(\ub_i)
\rb\TT(\ub_i)$. This strategy makes a compromise between
optimization (maximization of $f(\ub,\mthb^k)$, for $\ma_k$ small)
and estimation ($D$-optimal design, for $\ma_k$ large). Using the
results of Section \ref{S:results-LS}, the following is proved in
\cite{Pronzato00a} for LS estimation. When the errors $\mve_k$
form a martingale difference sequence with $\sup_k
\Ex\{\mve_k^2|\SF_{k-1}\}<\infty$, if $(\ma_k/k)\log \ma_k$ is
monotonically decreasing and $\ma_k/(\log k)^{1+\delta}$
monotonically increases to infinity for some $\delta>0$, then
$\mthb_{LS}^k \raas \mtbb$, $(1/k)\sum_{i=1}^k f(\ub_i,\mtbb)
\raas f(\ub^*,\mtbb)=\max_\ub f(\ub,\mtbb)$ and $\xi_k \raas
\delta_{\ub^*}$ (weak convergence of probability measures) as
$k\ra\infty$. That is, the LS estimator is strongly consistent,
and at the same time the design points $\ub_k$ tend to concentrate
at the optimal location $\ub^*$. Using Bayesian imbedding, the
same results are obtained when the conditions above on $\ma_k$ are
relaxed to $\ma_k\ra\infty$, $\ma_k/k\ra 0$, provided the errors
$\mve_k$ are i.i.d.\ $\SN(0,\ms^2)$, see \cite{PTa03}.

\subsection{Finite horizon: dynamic programming and dual
control}\label{S:finiteN}

The presentation is for self-tuning optimization, but the problem
is similar for other adaptive control situations. Suppose one
wishes to maximize $\sum_{i=1}^N w_i f(\ub_i,\mtbb)$ for some
sequence of positive weights $w_i$, with $\mtbb$ unknown and
estimated through observations $y_i=\eta(\mtbb,\ub_i)+\mve_i$. Let
$\pi$ denote a prior probability measure for $\bm\mt$ and define
$U_1^k=(\ub_1,\ldots,\ub_k)$, $Y_1^k=(y_1,\ldots,y_k)$ for all
$k$. The problem to be solved can then be written as
\begin{eqnarray}
&&  \max_{\ub_{1}}  \Ex\{w_1f(\ub_{1},\bm\mt) +
 \max_{\ub_2} \Ex\{w_2f(\ub_2,\bm\mt)  + \cdots \nonumber
\\
&&  \hspace{1cm} \max_{\ub_{N-1}} \Ex\{w_{N-1}f(\ub_{N-1},\bm\mt) \nonumber \\
&& \hspace{1cm} +  \max_{\ub_N}
\Ex\{w_Nf(\ub_N,\bm\mt)|U_1^{N-1},Y_1^{N-1} \}
\nonumber \\
&& \hspace{2cm} |U_1^{N-2},Y_1^{N-2} \} \cdots | \ub_1,y_1 \}  \}
 \label{SDP}
\end{eqnarray}
and thus corresponds to a Stochastic Dynamic Programming (SDP)
problem. It is, in general, extremely difficult to solve due to
the presence of imbedded maximizations and expectations. The
control $\ub_k$ has a {\em dual effect} (see e.g.
\cite{Bar-ShalomT74}): it affects both the value of
$f(\ub_k,\bm\mt)$ and the future uncertainty on $\bm\mt$ through
the posterior measures $\pi(\bm\mt|U_1^i,Y_1^i)$, $i\geq k$. One
of the main obstacles being the propagation of these measures,
classical approaches are based on their approximation. Consider
stage $k$, where $U_1^k$ and $Y_1^k$ are known. Then:\\
$\bullet$ Forced Certainty Equivalence control (FCE) replaces
$\pi(\bm\mt|U_1^i,Y_1^i)$
    for $i\geq k$ (a ``future posterior'' for $i>k$), by the delta measure
    $\delta_{\mth^k}$, where $\mthb^k$ is the current estimated
    value of $\bm\mt$ (see the examples of Sections \ref{S:STR} and
\ref{S:STO}); \\
$\bullet$ Open-Loop-Feedback-Optimal control (OLFO) replaces
$\pi(\bm\mt|U_1^i,Y_1^i)$, $i\geq k$,
    by the current posterior measure $\pi(\bm\mt|U_1^k,Y_1^k)$ (moreover, most often this
    posterior is approximated by a normal distribution $\SN(\mthb^k,
    \Cb_k)$).

The FCE and OLFO control strategies can be considered as {\em
passive} since they ignore the influence of
$\ub_{k+1},\ub_{k+2}\ldots$ on the future posteriors
$\pi(\bm\mt|U_1^i,Y_1^i)$, see, e.g., \cite{TseB-SM73}. On the
other hand, they yield a drastic simplification of the problem,
since the approximation of $\pi(\bm\mt|U_1^i,Y_1^i)$ for $i>k$
does not depend on the future observations $y_{k+1},y_{k+2}\ldots$
This, and the fact that few {\em active} alternatives exist,
explains their frequent usage.

The {\em active-control} strategy suggested in \cite{TseB-S73} is
based on a linearization around a nominal trajectory $\hat\ub(i)$
and extended Kalman filtering. It does not seem to have been much
employed, probably due to its rather high complexity. A
modification of OLFO control is proposed in \cite{PKWa96}. It
takes a very simple form when the model response
$\eta(\bm\mt,\ub)$ is linear in $\bm\mt$, that is,
$\eta(\bm\mt,\ub)=\rb\TT(\ub)\bm\mt$, the errors are i.i.d.\
normal $\SN(0,\ms^2)$ and the prior for $\bm\mt$ is also normal.
Then, at stage $k$, the posterior $\pi(\bm\mt|U_1^i,Y_1^i)$ is the
normal $\SN(\mthb^k_B, \Cb_k)$ for $i=k$ and can be approximated
by $\SN(\mthb^k_B, \Cb_i)$ for $i>k$, where $\mthb^k_B$ and
$\Cb_k$ are known (computed by classical recursive LS) and $\Cb_i$
follows a recursion similar to that of recursive LS,
$$
\Cb_{i+1}=\Cb_{i} -
\frac{\Cb_{i}\rb(\ub_{i+1})\rb\TT(\ub_{i+1})\Cb_{i}}{\ms^2+\rb\TT(\ub_{i+1})\Cb_{i}\rb(\ub_{i+1})}
\,, \ \ i\geq k \,.
$$
Note that $\Cb_i$ depends of $\ub_{k+1},\ub_{k+2}\ldots,\ub_i$
(which makes the strategy active), but not on
$y_{k+1},y_{k+2}\ldots$ (which makes it implementable). This
method has been successfully applied to the adaptive control of
model with a FIR, ARX, or state-space structure, see, e.g.,
\cite{KPWa96}, \cite{KPWc96}. It requires, however, that the
objective function $f(\ub,\bm\mt)$ in (\ref{SDP}) be non linear in
$\bm\mt$ to express the dependence in the covariance matrices
$\Cb_i$. Indeed, suppose that in the self-tuning optimization
problem the function to be maximized is the model response itself,
that is, $f(\ub,\bm\mt)=\rb\TT(\ub)\bm\mt$. Then,
$\Ex\{f(\ub,\bm\mt)|U_1^i,Y_1^i\}=\rb\TT(\ub)\mthb^i_B$ and using
the approximation $\SN(\mthb^k_B, \Cb_i)$ for the future
posteriors $\pi(\bm\mt|U_1^i,Y_1^i)$, $i>k$, one gets classical
FCE control based on the Bayesian estimator $\mthb^k_B$. On the
other hand, it is possible in that case to take benefit of the
linearity of the function and obtain an approximation of $\Ex\{
\max_\ub \rb\TT(\ub)\mthb^{N-1}_B | U_1^{N-2},Y_1^{N-2} \}$ for
small $\ms^2$, which can then be back-propagated; see \cite{PTa03}
where a control strategy is given that is within $\SO(\ms^4)$ of
the optimal (unknown) strategy $\ub_k^*$ for the SDP problem
(\ref{SDP}).

\section{Sequential DOE} \label{S:sequential-design}

Consider a nonlinear regression model for which the optimal design
problem consists in minimizing
$\Psi_{\mtb}(U_1^N)=\Phi[\Mb_F(U_1^N,\mtbb)]$ for some criterion
$\Phi$, with $\mtbb$ unknown. In order to design an experiment
adapted to $\mtbb$, a natural approach consists in working
sequentially. In {\em full-sequential design}, one support point
$\ub_k$ is introduced after each observation: $\mthb^{k-1}$ is
estimated from the data $(Y_1^{k-1},U_1^{k-1})$ and next $\ub_k$
minimizes $\Phi[\Mb_F(\{ U_1^{k-1},\ub_k\}, \mthb^{k-1})]$ (for
$D$-optimal design, this is equivalent to choosing $\ub_k$ that
maximizes $d_{\mth^{k-1}}(\ub_k,\xi_{k-1})$ with $d_\mt(\ub,\xi)$
the function (\ref{d}) and $\xi_{k-1}$ the empirical measure for
the design points in $U_1^{k-1}$). Note that it may be considered
as a FCE control strategy, where the input (design point) at step
$k$ is based on the current estimated value $\mthb^{k-1}$. For a
finite horizon $N$ (the number of observations), the problem is
similar to that of Section \ref{S:finiteN} (self-tuning
optimizer), with the design objective $\Phi[\Mb_F(U_1^N,\bm\mt)]$
substituted for $\sum_{i=1}^N w_if(\ub_i,\bm\mt)$. Although the
objective does not take an additive form, the problem is still of
the SDP type, and active-control strategies can thus be
constructed to approximate the optimal solution. However, they
seem to only provide marginal improvements over traditional
passive strategies like FCE control, see e.g.\
\cite{GPa00}\footnote{An active strategy aims at taking into
account the influence of current decisions on the future precision
of estimates; in that sense {\em DOE is naturally active by
definition}, even if based on FCE control. Trying to make
sequential DOE more active is thus doomed to small improvements}.

Although a sequentially designed experiment for the minimization
of $\Phi[\Mb_F(U_1^N,\bm\mt)]$ aims at estimating $\bm\mt$ with
maximum possible precision, it is difficult to assess that
$\mthb^N \raas \mtbb$ as $N\ra\infty$ (and thus $\xi_N \raas
\xi^*(\mtbb)$, with $\xi^*(\mtbb)$ the optimal design for $\mtbb$)
when full-sequential design is used; see \cite{Wu85} for a simple
example (with a positive answer) for LS estimation. When
full-sequential design is based on Bayesian estimation (posterior
mean), strong consistency can be proved if the optimal design
$\xi^*(\bm\mt)$ satisfies an identifiability condition for any
$\bm\mt$, see \cite{Hu98} (this is related to Bayesian imbedding
considered in Section \ref{S:results-bayesian}). The asymptotic
analysis of multi-stage sequential design is considered in
\cite{ChaudhuriM93} and the construction of asymptotically optimal
sequential design strategies in \cite{Spokoinyi92}, where it is
shown that {\em using two stages is enough}. Practical experience
tends to confirm the good performance of two-stage procedures,
see, e.g., \cite{BarenthinJH2005}.

\section{Concluding remarks and perspectives in DOE} \label{S:further-issues}

{\em Correlated errors.} Few results exist on DOE in the presence
of correlated observations and one can refer e.g.\ to
\cite{PazmanM2001}, \cite{MullerP2003}, \cite{FedorovM2004} and
the monograph \cite{Muller2001} for recent developments. See also
Section \ref{S:Model-based}.
 The situation is different in the adaptive
control community where correlated errors are classical, see
Section \ref{S:results-LS} (for instance, the paper
\cite{Ninness2000} gives results on strong laws of large numbers
for correlated sequences of random variables under rather common
assumptions in signal or control applications), which calls for
appropriate developments in DOE. Notice that when the correlation
of the error process decays at hyperbolic rate (long-range
dependence), the asymptotic theory of parameter estimation in
regression models (Section \ref{S:information-matrix}) must itself
be revisited, see, e.g., \cite{IvanovL2004}.

{\em Nonlinear models.} The presentation in Sections
\ref{S:control-for-design} and \ref{S:adaptive-control} has
concerned models with a linear dynamic (e.g.\ with a Box and
Jenkins structure), but models that corresponds to nonlinear
differential or recurrence equations raise no special difficulty
for the construction of the Fisher information matrix (Section
\ref{S:control-for-design}), which can always be obtained through
simulations. The main issue concerns {\em linearity with respect
to the model parameters} $\bm\mt$. In particular, few results
exist concerning the extension of the results of Section
\ref{S:some-results} to models with a nonlinear parametrization
(see \cite{Lai94} for LS estimation and \cite{Hu98b} for results
on Bayesian imbedding when $\bm\mt$ has a discrete prior).

{\em DOE without persistence of excitation.} In the context of
self-tuning regulation, we mentioned in Section
\ref{S:results-bayesian} that random perturbations may be added to
certainty equivalence control based on LS estimation to guarantee
the strong consistency of parameter estimates and the
asymptotically optimal growth of the control objective, see
\cite{LaiW87}. This is an example of a situation where
``non-stationary experiments'' could be designed in order to
replace random perturbations by inputs with suitable spectrum and
asymptotically vanishing amplitude. In the same vein, the modified
OLFO control proposed in \cite{PKWa96} and the small-noise
approximation of \cite{PTa03} (designed for self-tuning
optimization, but extendable to self-tuning regulation) make a
good compromise between exploration and exploitation when the
horizon is finite (see Section \ref{S:finiteN}). An asymptotic
analysis for the horizon tending to infinity could permit to
design simpler non-stationary strategies.

{\em Nonparametric models, active learning and control.}
Strategies are called active in opposition to passive ones that
collect data ``as they come''. DOE is thus intrinsically active,
and its use in learning leads to methods that try to select
training samples instead of taking them randomly. Although its
usefulness is now well perceived in the statistical learning
community, it is still at an early stage of development due to the
complexity of DOE for nonparametric models. More generally, active
strategies are valuable each time actions or decisions have a dual
effect and a compromise should be made between exploration and
exploitation: exploration may be done at random, but better
performance is achieved when it is carefully planned. For
instance, active strategies connected with Markov decision theory
could yield improvements in reinforcement learning, see e.g. the
survey \cite{KaelblingLM96}.

Linking nonparametric estimation with control forms a quite
challenging area, where the issues raised by the estimation of the
function that defines the dynamics of the system come in addition
to those, more classical, of adaptive control with parametric
models, see, e.g., \cite{PortierO2000}, \cite{PoggiP2000} for
emerging developments.

{\em Algorithms for optimal DOE.} The importance of constructing
criteria for DOE in relation with the intended objective has been
underlined in Section \ref{S:frequency-domain} where criteria of
the minimax type have been introduced from robust-control
considerations. (Minimax-optimal design is also an efficient
method to face the dependence of local optimal design in the
unknown value of the model parameters, see Section
\ref{S:design-nonlinear}.) Although the minimax problem can often
be formulated as a convex one, sometimes with a finite number of
constraints, the development of specific algorithms would be much
profitable to the diffusion of the methodology, in the same way as
the classical design algorithms of Sections
\ref{S:algorithms-exact} and \ref{S:algorithms-for-approximateD}
have contributed to the diffusion of optimal DOE outside the
statistical community where it originated.

{\em Another view on global optimization.} Let $f(\ub)$ be a
function to be maximized with respect to $\ub$ in some given set
$\SU$; it is not assumed to be concave, nor is the set $\SU$
assumed to be convex, so that local maxima may exist. The function
can be evaluated at any given input $\ub_i\in\SU$, which gives an
``observation'' $y(\ub_i)=f(\ub_i)$. In engineering applications
where the evaluation of $f$ corresponds to the execution of a
large simulation code, expensive in terms of computing time, it is
of paramount importance to use an optimization method parsimonious
in terms of number of function evaluations. This enters into the
framework of {\em computer experiments}, where Kriging is now a
rather well-established tool for modelling, see Section
\ref{S:Kriging}. Using a Bayesian point of view, the value
$f(\ub)$ after the collection of the data
$\SD_k=\{[\ub_1,y(\ub_1)],\ldots,[\ub_k,y(\ub_k)]\}$ can be
considered as distributed with the density $\varphi(y|\SD_k,\ub)$
of the normal distribution
$\SN(\yh_{\SD_k}(\ub),\rho^2_{\SD_k}(\ub))$, where
$\yh_{\SD}(\ub)$ and $\rho^2_{\SD}(\ub)$ are respectively given by
(\ref{yh-kriging}) and (\ref{MSE-kriging}). An optimization
algorithm that uses this information should then make a compromise
between exploration (trying to reduce the MSE $\rho^2_{\SD}(\ub)$
by placing observations at values of $\ub$ where
$\rho^2_{\SD}(\ub)$ is large) and exploitation (trying to maximize
the expected response $\yh_{\SD}(\ub)$). A rather intuitive method
is to choose $\ub_{k+1}$ that maximizes $\yh_{\SD_k}(\ub)+ \alpha
\rho_{\SD_k}(\ub)$ for some positive constant $\alpha$, see
\cite{CoxJ92}. In theory, Stochastic Dynamic Programming could be
used to find the optimal strategy (or algorithm) to maximize
$f(\ub)$: when the number $N$ of evaluations is given in advance,
the problem takes the same form as in (\ref{SDP}) with $w_i=0$ for
$i=1,\ldots,N-1$. However, in practise this SDP problem is much
too difficult to solve, and approximations must be used to define
suboptimal searching rules. For instance, one may use a
one-step-ahead approach and choose the input $\ub_{k+1}$ that
maximizes the {\em expected improvement} $EI(\ub) =
\int^{\infty}_{y_k^{\max}} [y-y_k^{\max}] \, \varphi(y|\SD_k,\ub)
\, dy$, with $y_k^{\max}=\max\{y(\ub_1),\ldots,y(\ub_k)\}$, the
maximum value of $f$ observed so far, see \cite{MockusTZ78},
\cite{Mockus89}, \cite{SchonlauWJ98}. The function $f$ is then
evaluated at $\ub_{k+1}$, the Kriging model is updated (although
not necessarily at each iteration), and similar steps are
repeated. Each iteration of the resulting algorithm requires one
evaluation of $f$ and the solution of another global optimization
problem, for which any ad-hoc global search algorithm can be used
(the optimization concerns the function $EI$, which is easier to
evaluate than $f$). For instance, it is suggested in \cite{BPc01}
to update a Delaunay triangulation of the set $\SU$ based on the
vertices $\ub_i$, $i=1,\ldots,k$, and to perform the global search
for the maximization of $EI(\ub)$ by initializing local searches
at the centers of the Delaunay cells. Note that the algorithm
tends asymptotically to observe everywhere in $\SU$, since the
expected improvement $EI(\ub)$ is always strictly positive at any
value $\ub$ where no observation has been taken yet. However, a
credible stopping rule is given by the criterion itself: it is
reasonable to stop when the expected improvement becomes small
enough. One can refer to \cite{SchonlauWJ98}, \cite{JonesSW98} for
detailed implementations, including problems with constraints also
defined by simulation codes. Derivative information on $f$ can be
included in the Kriging model, as indicated in Section
\ref{S:Kriging}, and thus used by the optimization algorithm, see
\cite{LearyBK2004}. It seems that suboptimal searching rules
looking further than one-step-ahead have never been used, which
forms a rather challenging objective for active control. Also, the
one-step-ahead method above is completely passive with respect to
the estimation of the parameters of the Kriging model, and active
strategies (even one-step-ahead) are still to be designed, see
Section \ref{S:Model-based}. Note finally that the definition of
the expected improvement $EI$ is not adapted to the presence of
errors in the evaluation of $f$, so that further developments are
required for situations where one optimizes the observed response
of a real physical process.

\begin{figure}[htb]
\begin{center}
\includegraphics[scale=0.35]{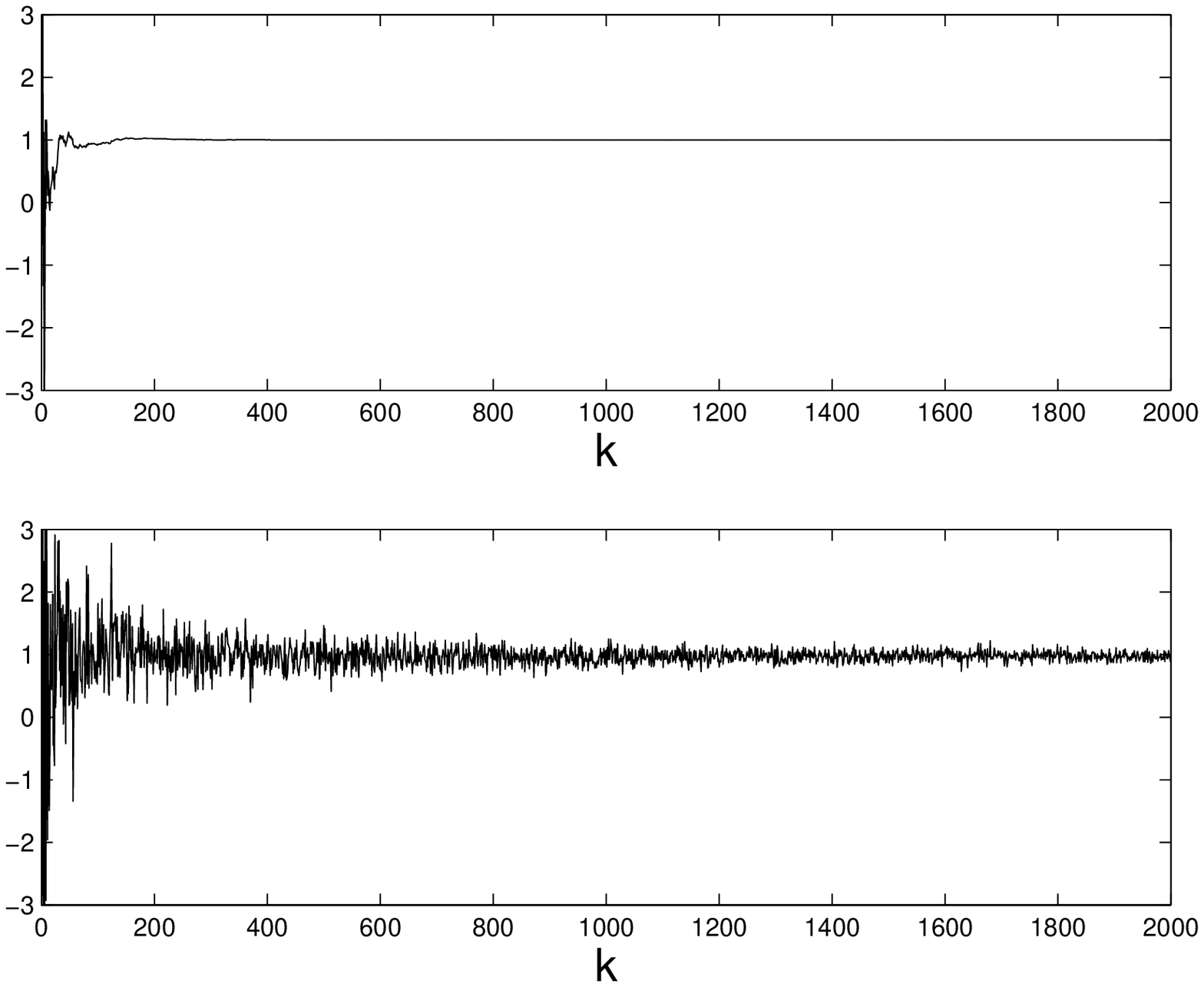}
\caption{Top: evolution of the LS estimator $\mth^k$; bottom:
evolution of $\tilde\mt^k$ given by (\ref{tildemt}).}
\label{F:KKK-thetas}
\end{center}
\end{figure}

{\em NFC, FCE, estimating functions and DOE.} Consider again the
example of NFC in Section \ref{S:NFC}, in the case where the state
$x_k$ is observed though $y_k=x_k+\mve_k$ and $y_k$ is substituted
for $x_k$ in (\ref{NFC-nonoise}). As shown in Figure
\ref{F:KKK-nonoise}, the NFC estimator $\mth_k$ is then not
consistent. On the other hand, in the same situation more
classical estimation techniques have satisfactory behaviors
(hence, in terms of DOE, NFC brings enough information to estimate
the unknown $\mtb$). Consider for instance the LS estimator of
$\mtb$ in (\ref{discrete}), obtained from $y_1,\ldots,y_k$. Since
$x_k$ is nonlinear in $\mt$, recursive LS cannot be used
directly\label{PFN:processnoise}\footnote{\label{FN:processnoise}
The situation would be much easier for the autoregressive model
$y_{i+1}=y_i + T[ u_i + \mtb(y_i+1)]+ \mve_{i+1}$ where $\mt$
could be estimated by recursive LS. Note that this model can be
considered as resulting from the discretization of $\dot y = u +
\mtb(y+1) + S dB_t(\omega)$, with $B_t(\omega)$ the standard
Brownian motion (starting at zero with variance 1), which
corresponds to the introduction of process noise into
(\ref{diff-eq}), and gives $\mve_{i+1}= S
[B_{(i+1)T}(\omega)-B_{iT}(\omega)]$.}, but the estimation becomes
almost recursive using a stochastic-Newton algorithm, see, e.g.,
\cite{WPl97}, p.\ 208. Figure \ref{F:KKK-thetas} (top) shows that
the corresponding estimator $\mth^k$ converges quickly to the true
value $\mtb=1$. The evolution of the estimator (\ref{tildemt})
obtained from an estimating function is presented on the same
figure (bottom). Its convergence is slower than that of $\mth^k$,
due in particular to the presence of the term $y_k/(kT)$ in the
numerator of (\ref{tildemt}), but {\em its construction is much
simpler}. An important consequence is that the analysis of its
asymptotic behavior in an adaptive-control framework is easier
than for LS estimation: $\tilde\mt^k$ is a consistent estimator of
$\mtb$ when $(1/k) \sum_{i=1}^{k-1} x_i$ is bounded away from $-1$
(which is the case since the control drives this quantity to
zero). Simulations confirm that when applying the FCE controller
$u_k=-(a+ \tilde\mt^k) \hat x_k(\tilde\mt^k)- \tilde\mt^k$ to
(\ref{discrete}), with $\hat x_k(\tilde\mt^k)$ obtained by
substituting $\tilde\mt^k$ for $\mtb$ in (\ref{discrete}),
$\tilde\mt^k$ converges to $\mtb$ and the state $x_k$ is correctly
driven to zero. Simulations show, however, that the dynamic of the
state is slower than for NFC; it is thus tempting to combine both
strategies. For instance, one could use FCE control based on
$\tilde\mt^k$ when the standard deviation of $\tilde\mt^k$ is
smaller than some prescribed value, and use NFC otherwise. The
simulation results obtained with such a switching strategy are
encouraging and indicate that the combination of different
estimation methods may improve the performance of the controller.
At the same time, this raises more issues than it brings answers.
Some are listed below.

Combining NFC, which relies on Lyapunov stability, with simple
predictors, e.g. based on FCE, while preserving stability, forms
an interesting challenge for which results on input-to-state
stabilizing control could be used (see, e.g., Chapters 5 and 6 of
\cite{KrsticKK95} for continuous-time and \cite{NesicL2002} for
discrete-time control). In classical FCE control the consistency
of the estimator is a major issue. Using suitable estimating
functions could then lead to fruitful developments, due the
flexibility of the method and the simplicity of the associated
estimators. Suitably designed perturbations could be introduced
for helping the estimation, possibly following developments
similar to those that lead to the active-control strategies of
Sections \ref{S:results-bayesian} and \ref{S:finiteN}. At the same
time, the perturbed control should not endanger the stability of
the system. Designing input sequences (possibly vanishing with
time) that bring maximum information for estimation subject to a
stability constraint forms an unusual and challenging DOE problem.
Finally, as a first step towards the design of robust-and-adaptive
controllers mentioned at the end of Section
\ref{S:frequency-domain}, one may replace a traditional FCE
controller by one that gives the best performance for the worst
model in the current uncertainty set (roughly speaking, for the
self-tuning problem considered in Section \ref{S:finiteN} this
amounts to replacing expectations in (\ref{SDP}) by minimizations
with respect to $\bm\mt$ in the current uncertainty set). The
determination of active-control strategies for such minimax
(dynamical games) problems seems to be a promising direction for
developments in adaptive control.

\begin{ack}
This work was partially supported by the IST Programme of the
European Community, under the PASCAL Network of Excellence,
IST-2002-506778. This publication only reflects the authors's
view. The paper is partly based on a three-hours mini-course given
at the 24th Benelux Meeting at Houffalize, Belgium, on March
22-24, 2005. It has benefited from several interesting and
motivating discussions during that meeting, and the organizers of
the event are gratefully acknowledged for their invitation. Many
ideas and references result from many years of collaboration with
\'Eric Walter (CNRS/SUPELEC/Universit\'e Paris XI, France), Andrej
P\'azman (Comenius University, Bratislava, Slovakia), Henry P.\
Wynn (LSE, London, UK) and Anatoly A. Zhigljavsky (Cardiff
University, UK). Sections \ref{S:results-bayesian} and
\ref{S:further-issues} have benefited from several discussions
with \'Eric Thierry and Tarek Hamel at I3S. The encouraging
comments and suggestions of several referees were also much
helpful.
\end{ack}


\bibliographystyle{plain}        

\begin{thebibliography}{100}

\bibitem{AstromW73}
K.J. {{\AA}str\"{o}m} and B.~Wittenmark.
\newblock On self-tuning regulators.
\newblock {\em Automatica}, 9:195--199, 1973.

\bibitem{AstromW89}
K.J. {{\AA}str\"{o}m} and B.~Wittenmark.
\newblock {\em Adaptive Control}.
\newblock Addison Wesley, 1989.

\bibitem{AtkinsonC74}
A.C. Atkinson and D.R. Cox.
\newblock Planning experiments for discriminating between models (with
  discussion).
\newblock {\em Journal of Royal Statistical Society}, B36:321--348, 1974.

\bibitem{AtkinsonD92}
A.C. Atkinson and A.N. Donev.
\newblock {\em Optimum Experimental Design}.
\newblock Oxford University Press, 1992.

\bibitem{AtkinsonF75a}
A.C. Atkinson and V.V. Fedorov.
\newblock The design of experiments for discriminating between two rival
  models.
\newblock {\em Biometrika}, 62(1):57--70, 1975.

\bibitem{AtkinsonF75b}
A.C. Atkinson and V.V. Fedorov.
\newblock Optimal design: experiments for discriminating between several
  models.
\newblock {\em Biometrika}, 62(2):289--303, 1975.

\bibitem{Bar-ShalomT74}
Y.~Bar-Shalom and E.~Tse.
\newblock Dual effect, certainty equivalence, and separation in stochastic
  control.
\newblock {\em IEEE Transactions on Automatic Control}, 19(5):494--500, 1974.

\bibitem{BarenthinJH2005}
M.~Barenthin, H.~Jansson, and H.~Hjalmarsson.
\newblock Applications of mixed ${H}_2$ and ${H}_\infty$ input design in
  identification.
\newblock In {\em Proc. 16th IFAC World Congress on Automatic Control}, Prague,
  2005.
\newblock CD-ROM -- Paper 03882.

\bibitem{BarkerG99}
H.A. Barker and K.R. Godfrey.
\newblock System identification with multi-level periodic perturbation signals.
\newblock {\em Control Engineering Practice}, 7:717--726, 1999.

\bibitem{Bartlett2003}
P.L. Bartlett.
\newblock Prediction algorithms: complexity, concentration and convexity.
\newblock In {\em Prep.\ 13th IFAC Symposium on System Identification,
  Rotterdam}, pages 1507--1517, August 2003.

\bibitem{BPc01}
R.~Bates and L.~Pronzato.
\newblock Emulator-based global optimisation using lattices and {D}elaunay
  tesselation.
\newblock In P.~Prado and R.~Bolado, editors, {\em Proc.\ 3rd Int.\ Symp.\ on
  Sensitivity Analysis of Model Output}, pages 189--192, Madrid, June 2001.

\bibitem{BiedermannD2001}
S.~Biedermann and H.~Dette.
\newblock Minimax optimal designs for nonparametric regression --- a further
  optimality property of the uniform distribution.
\newblock In P.~Hackl A.C.~Atkinson and W.G. M{\"u}ller, editors, {\em mODa'6
  -- Advances in Model--Oriented Design and Analysis, Proceedings of the 76th
  Int.\ Workshop, Puchberg/Schneberg (Austria)}, pages 13--20, Heidelberg, June
  2001. Physica Verlag.

\bibitem{Bohning89}
D.~B{\"o}hning.
\newblock Likelihood inference for mixtures: geometrical and other
  constructions of monotone step-length algorithms.
\newblock {\em Biometrika}, 76(2):375--383, 1989.

\bibitem{BomboisAG2005}
X.~Bombois,~B.D.O.~Anderson,~and~M.~Gevers.
Quantification of frequency domain error bounds with guaranteed
  confidence level in prediction error identification.
\newblock {\em System \& Control Letters}, 54:471--482, 2005.

\bibitem{BomboisSGHVdH2004}
X.~Bombois, G.~Scorletti, M.~Gevers, R.~Hildebrand, and P.~{Van den Hof}.
\newblock Cheapest open-loop identification for control.
\newblock In {\em Proc. 43rd Conf. on Decision and Control}, pages 382--387,
  The Bahamas, December 2004.

\bibitem{BomboisSGVdHH2005}
X.~Bombois, G.~Scorletti, M.~Gevers, P.~{Van den Hof}, and R.~Hildebrand.
\newblock Least costly identification experiment for control.
\newblock {\em Automatica}, 42(10):1651--1662, 2006.

\bibitem{BoxH67}
G.E.P. Box and W.J. Hill.
\newblock Discrimination among mechanistic models.
\newblock {\em Technometrics}, 9(1):57--71, 1967.

\bibitem{BoxW51}
G.E.P. Box and K.B. Wilson.
\newblock On the experimental attainment of optimum conditions (with
  discussion).
\newblock {\em Journal of Royal Statistical Society}, B13(1):1--45, 1951.

\bibitem{Box71}
M.J. Box.
\newblock Bias in nonlinear estimation.
\newblock {\em Journal of Royal Statistical Society}, B33:171--201, 1971.

\bibitem{BoydEGFB94}
S.~Boyd, L.~{E}l {G}haoui, E.~Feron, and V.~Balakrishnan.
\newblock {\em Linear Matrix Inequalities in System and Control Theory}.
\newblock SIAM, Philadelphia, 1994.

\bibitem{BoydV2004}
S.~Boyd and L.~Vandenberghe.
\newblock {\em Convex Optimization}.
\newblock Cambridge University Press, Cambridge, 2004.

\bibitem{BozinZ91}
A.S. Bozin and M.B. Zarrop.
\newblock Self tuning optimizer --- convergence and robustness properties.
\newblock In {\em Proc. 1st European Control Conf.}, pages 672--677, Grenoble,
  July 1991.

\bibitem{Breiman96}
L.~Breiman.
\newblock Bagging predictors.
\newblock {\em Machine Learning}, 24(2):123--140, 1996.

\bibitem{Caines88}
P.E. Caines.
\newblock {\em Linear Stochastic Systems}.
\newblock Wiley, New York, 1988.

\bibitem{CampiW2005}
M.C. Campi and E.~Weyer.
\newblock Guaranteed non-asymptotic confidence regions in system
  identification.
\newblock {\em Automatica}, 41(10):1751--1764, 2005.

\bibitem{ChalonerL89}
K.~Chaloner and K.~Larntz.
\newblock Optimal {B}ayesian design applied to logistic regression experiments.
\newblock {\em Journal of Statistical Planning and Inference}, 21:191--208,
  1989.

\bibitem{ChalonerV95}
K.~Chaloner and I.~Verdinelli.
\newblock {B}ayesian experimental design: a review.
\newblock {\em Statistical Science}, 10(3):273--304, 1995.

\bibitem{ChaudhuriM93}
P.~Chaudhuri and P.A. Mykland.
\newblock Nonlinear experiments: optimal design and inference based likelihood.
\newblock {\em Journal of the American Statistical Association},
  88(422):538--546, 1993.

\bibitem{Chen80}
C.S. Chen.
\newblock Optimality of some weighing and $2^n$ fractional factorial designs.
\newblock {\em Annals of Statistics}, 8:436--446, 1980.

\bibitem{ChengHT98}
M.-Y Cheng, P.~Hall, and D.M. Titterington.
\newblock Optimal design for curve estimation by local linear smoothing.
\newblock {\em Bernoulli}, 4(1):3--14, 1998.

\bibitem{Chernoff53}
H.~Chernoff.
\newblock Locally optimal designs for estimating parameters.
\newblock {\em Annals of Math. Stat.}, 24:586--602, 1953.

\bibitem{ChoiKAL2002}
J.-Y. Choi, M.~Krsti\'c, K.B. Ariyur, and J.S. Lee.
\newblock Extremum seeking control for discrete-time systems.
\newblock {\em IEEE Transactions on Automatic Control}, 47(2):318--323, 2002.

\bibitem{Cohn94}
D.A. Cohn.
\newblock Neural network exploration using optimal experiment design.
\newblock In J.~Cowan, G.~Tesauro, and J.~Alspector, editors, {\em Advances in
  Neural Information Processing Systems 6}. Morgan Kaufmann, 1994.

\bibitem{CohnGJ96}
D.A. Cohn, Z.~Ghahramani, and M.I. Jordan.
\newblock Active learning with statistical models.
\newblock {\em Journal of Artificial Intelligence Research}, 4:129--145, 1996.

\bibitem{CoxJ92}
D.D. Cox and S.~John.
\newblock A statistical method for global optimization.
\newblock In {\em Proc. IEEE Int. Conf. on Systems Man and Cybernetics},
  volume~2, pages 1241--1246, Chicago, IL, October 1992.

\bibitem{CoxH74}
D.R. Cox and D.V. Hinkley.
\newblock {\em Theoretical Statistics}.
\newblock Chapman \& Hall, London, 1974.

\bibitem{CuckerS2001}
F.~Cucker and S.~Smale.
\newblock On the mathematical foundations of learning.
\newblock {\em Bulletin of the {AMS}}, 39(1):1--49, 2001.

\bibitem{CurrinMMY91}
C.~Currin, T.J. Mitchell, M.D. Morris, and D.~Ylvisaker.
\newblock {B}ayesian prediction of deterministic functions, with applications
  to the design and analysis of computer experiments.
\newblock {\em Journal of the American Statistical Association}, 86:953--963,
  1991.

\bibitem{D'Argenio81}
D.Z. D'Argenio.
\newblock Optimal sampling times for pharmacokinetic experiments.
\newblock {\em Journal of Pharmacokinetics and Biopharmaceutics},
  9(6):739--756, 1981.

\bibitem{denHertog94}
D.~{den Hertog}.
\newblock {\em Interior Point Approach to Linear, Quadratic and Convex
  Programming}.
\newblock Kluwer, Dordrecht, 1994.

\bibitem{Faraway90}
J.J. Faraway.
\newblock Sequential design for the nonparametric regression of curves and
  surfaces.
\newblock {\em Computing Science and Statistics}, 22:104--110, 1990.

\bibitem{Fedorov72}
V.V. Fedorov.
\newblock {\em Theory of Optimal Experiments}.
\newblock Academic Press, New York, 1972.

\bibitem{Fedorov80}
V.V. Fedorov.
\newblock Convex design theory.
\newblock {\em Math. Operationsforsch. Statist., Ser. Statistics},
  11(3):403--413, 1980.

\bibitem{FedorovH97}
V.V. Fedorov and P.~Hackl.
\newblock {\em Model-Oriented Design of Experiments}.
\newblock Springer, Berlin, 1997.

\bibitem{FedorovM2004}
V.V.~Fedorov~and~W.G.~M{\"u}ller.~Optimum~design~for correlated~processes~via~eigenfunction~expansions.
Department of Statistics and Mathematics, Wirt\-schafts\-universit{\"a}t
  Wien, Report 6, June 2004.

\bibitem{Fisher25}
R.A. Fisher.
\newblock {\em Statistical Methods for Research Workers}.
\newblock Oliver \& Boyd, Edimbourgh, 1925.

\bibitem{FliessS-R2003}
M.~Fliess and H.~Sira-{R}ami\'irez.
\newblock An algebraic framework for linear identification.
\newblock {\em ESAIM: Control, Optimization and Calculus of Variations},
  (9):151--168, 2003.

\bibitem{ForsgrenGW2002}
A.~Forsgren, P.E. Gill, and M.H. Wright.
\newblock Interior methods for nonlinear optimization.
\newblock {\em SIAM Review}, 44(4):525--597, 2002.

\bibitem{ForssellL99}
U.~Forssell and L.~Ljung.
\newblock Closed-loop identification revisited.
\newblock {\em Automatica}, 35:1215--1241, 1999.

\bibitem{ForssellL2000}
U.~Forssell and L.~Ljung.
\newblock Some results on optimum experiment design.
\newblock {\em Automatica}, 36:749--756, 2000.

\bibitem{GPa00}
R.~Gautier and L.~Pronzato.
\newblock Adaptive control for sequential design.
\newblock {\em Discussiones Mathematicae, Probability \& Statistics},
  20(1):97--114, 2000.

\bibitem{GazutMI2006}
S.~Gazut, J.-M. Martinez, and S.~Issanchou.
\newblock Plans d'exp\'eriences it\'eratifs pour la construction de mod\`eles
  non lin\'eaires.
\newblock In {\em CD -- 38\`emes Journ\'ees de Statistique de la SFdS},
  Clamart, France, 2006.

\bibitem{Gevers2005}
M.~Gevers.
\newblock Identification for control. {F}rom the early achievements to the
  revival of experimental design.
\newblock {\em European Journal of Control}, 11(45):335--352, 2005.

\bibitem{GeversBCSA2003a}
M.~Gevers, X.~Bombois, B.~Codrons, G.~Scorletti, and B.D.O. Anderson.
\newblock Model validation for control and controller validation in a
  prediction error identification framework --- {P}art {I}: theory.
\newblock {\em Automatica}, 39:403--415, 2003.

\bibitem{GeversBCSA2003b}
M.~Gevers, X.~Bombois, B.~Codrons, G.~Scorletti, and B.D.O. Anderson.
\newblock Model validation for control and controller validation in a
  prediction error identification framework --- {P}art {II}: illustrations.
\newblock {\em Automatica}, 39:417--427, 2003.

\bibitem{GeversL85}
M.~Gevers and L.~Ljung.
\newblock Benefits of feedback in experiment design.
\newblock In {\em Prep. 7th IFAC/IFORS Symp. on Identification and System
  Parameter Estimation}, pages 909--914, York, 1985.

\bibitem{GeversL86}
M.~Gevers and L.~Ljung.
\newblock Optimal experiment design with respect to the intended model
  application.
\newblock {\em Automatica}, 22:543--554, 1986.

\bibitem{GoodwinP77}
G.C. Goodwin and R.L. Payne.
\newblock {\em Dynamic System Identification: Experiment Design and Data
  Analysis}.
\newblock Academic Press, New York, 1977.

\bibitem{GoodwinS89}
G.C. Goodwin and M.E. Salgado.
\newblock A stochastic imbedding approach for quantifying uncertainty in the
  estimation of restricted complexity models.
\newblock {\em International Journal of Adaptive Control and Signal
  Processing}, 3:333--356, 1989.

\bibitem{Guo94}
L.~Guo.
\newblock Further results on least-squares based adaptive minimum variance
  control.
\newblock {\em SIAM J. Control and Optimization}, 32(1):187--212, 1994.

\bibitem{HPa06_SPL}
R.~Harman and L.~Pronzato.
\newblock Improvements on removing non-optimal support points in {D}-optimum
  design algorithms.
\newblock {\em Statistics \& Probability Letters}, 77:90--94, 2007.

\bibitem{HastieTF2001}
T.~Hastie, R.~Tibshirani, and J.~Friedman.
\newblock {\em The Elements of Statistical Learning. Data Mining, Inference and
  Prediction}.
\newblock Springer, Heidelberg, 2001.

\bibitem{Heyde97}
C.C. Heyde.
\newblock {\em Quasi-likelihood and its Application. A General Approach to
  Optimal Parameter Estimation}.
\newblock Springer, New York, 1997.

\bibitem{HildebrandG2003}
R.~Hildebrand and M.~Gevers.
\newblock Identification for control: optimal input design with respect to a
  worst-case $\nu$-gap cost function.
\newblock {\em SIAM Journal on Control and Optimization}, 42(5):1586--1608,
  2003.

\bibitem{Hill78}
P.D.H. Hill.
\newblock A review of experimental design procedures for regression model
  discrimination.
\newblock {\em Technometrics}, 20:15--21, 1978.

\bibitem{Hjalmarsson2005}
H.~Hjalmarsson.
\newblock From experiment design to closed-loop control.
\newblock {\em Automatica}, 41:393--438, 2005.

\bibitem{HjalmarssonGD-B96}
H.~Hjalmarsson, M.~Gevers, and F.~{D}e {B}ruyne.
\newblock For model-based control design, closed-loop identification gives
  better performance.
\newblock {\em Automatica}, 32(12):1659--1673, 1996.

\bibitem{Hu96}
I.~Hu.
\newblock Strong consistency of {B}ayes estimates in stochastic regression
  models.
\newblock {\em Journal of Multivariate Analysis}, 57:215--227, 1996.

\bibitem{Hu97}
I.~Hu.
\newblock Strong consistency in stochastic regression models via posterior
  covariance matrices.
\newblock {\em Biometrika}, 84(3):744--749, 1997.

\bibitem{Hu98}
I.~Hu.
\newblock On sequential designs in nonlinear problems.
\newblock {\em Biometrika}, 85(2):496--503, 1998.

\bibitem{Hu98b}
I.~Hu.
\newblock Strong consistency of {B}ayes estimates in nonlinear stochastic
  regression models.
\newblock {\em Journal of Statistical Planning and Inference}, 67:155--163,
  1998.

\bibitem{Huber81}
P.J. Huber.
\newblock {\em Robust Statistics}.
\newblock John Wiley, New York, 1981.

\bibitem{IvanovL2004}
A.V. Ivanov and N.N. Leonenko.
\newblock Asymptotic theory of nonlinear regression with long-range dependence.
\newblock {\em Mathematical Methods of Statistics}, 13(2):153--178, 2004.

\bibitem{JanssonH2004}
H.~Jansson and H.~Hjalmarsson.
\newblock Mixed ${H}_\infty$ and ${H}_2$ input design for identification.
\newblock In {\em Proc. 43rd Conf. on Decision and Control}, pages 388--393,
  The Bahamas, December 2004.

\bibitem{JanssonH2005}
H.~Jansson and H.~Hjalmarsson.
\newblock Input design via {LMI}s admitting frequency-wise model specifications
  in confidence regions.
\newblock {\em IEEE Transactions on Automatic Control}, 50(10):1534--1549,
  2005.

\bibitem{JanssonH-IFAC2005}
H.~Jansson and H.~Hjalmarsson.
\newblock Optimal experiment design in closed loop.
\newblock In {\em Proc. 16th IFAC World Congress on Automatic Control}, Prague,
  2005.
\newblock CD-ROM -- Paper 04528.

\bibitem{Jennrich69}
R.I. Jennrich.
\newblock Asymptotic properties of nonlinear least squares estimation.
\newblock {\em Annals of Math. Stat.}, 40:633--643, 1969.

\bibitem{JohnsonMY90}
M.E. Johnson, L.M. Moore, and D.~Ylvisaker.
\newblock Minimax and maximin distance designs.
\newblock {\em Journal of Statistical Planning and Inference}, 26:131--148,
  1990.

\bibitem{JonesSW98}
D.~Jones, M.~Schonlau, and W.J. Welch.
\newblock Efficient global optimization of expensive black-box functions.
\newblock {\em Journal of Global Optimization}, 13:455--492, 1998.

\bibitem{KaelblingLM96}
L.P. Kaelbling, M.L. Littman, and A.W. Moore.
\newblock Reinforcement learning: a survey.
\newblock {\em Journal of Artificial Intelligence Research}, 4:237--285, 1996.

\bibitem{Kanamori2002}
T.~Kanamori.
\newblock Statistical asymptotic theory of active learning.
\newblock {\em Annals Inst. Statist. Math.}, 54(3):459--475, 2002.

\bibitem{Kerestecioglu93}
F.~Kerestecio\u{g}lu.
\newblock {\em Change Detection and Input Design in Dynamical Systems}.
\newblock Research Studies Press, Taunton, UK, 1993.

\bibitem{KieferW52}
J.~Kiefer and J.~Wolfowitz.
\newblock Stochastic estimation of the maximum of a regression function.
\newblock {\em Annals of Math. Stat.}, 23:462--466, 1952.

\bibitem{KieferW59}
J.~Kiefer and J.~Wolfowitz.
\newblock Optimum designs in regression problems.
\newblock {\em Annals of Math. Stat.}, 30:271--294, 1959.

\bibitem{KieferW60}
J.~Kiefer and J.~Wolfowitz.
\newblock The equivalence of two extremum problems.
\newblock {\em Canadian Journal of Mathematics}, 12:363--366, 1960.

\bibitem{Krige51}
D.G. Krige.
\newblock A statistical approach to some mine valuation and allied problems on
  the {W}itwatersrand.
\newblock Master Thesis, University of Witwatersrand, 1951.

\bibitem{Krstic2000}
M.~Krsti\'c.
\newblock Performance improvement and limitations in extrremum seeking control.
\newblock {\em System \& Control Letters}, 39:313--326, 2000.

\bibitem{KrsticKK95}
M.~Krsti\'c, I.~Kanellakopoulos, and P.~Kokotovi\'c.
\newblock {\em Nonlinear and Adaptive Control Design}.
\newblock Wiley, New York, 1995.

\bibitem{KrsticW2000}
M.~Krsti\'c and H.-H. Wang.
\newblock Stability of extremum seeking feedback for general nonlinear dynamic
  systems.
\newblock {\em Automatica}, 36:595--601, 2000.

\bibitem{KubruslyM85}
C.S. Kubrusly and H.~Malebranche.
\newblock Sensors and controllers location in distributed systems---{A} survey.
\newblock {\em Automatica}, 21(2):117--128, 1985.

\bibitem{KPWa96}
C.~Kulcs\'ar, L.~Pronzato, and E.~Walter.
\newblock Dual control of linearly parameterised models via prediction of
  posterior densities.
\newblock {\em European J. of Control}, 2(1):135--143, 1996.

\bibitem{KPWc96}
C.~Kulcs\'ar, L.~Pronzato, and E.~Walter.
\newblock Experimental design for the control of linear state-space systems.
\newblock In {\em Proc. 13th IFAC World Congress}, volume~C, pages 175--180,
  San Francisco, June 1996.

\bibitem{Kumar90}
P.R. Kumar.
\newblock Convergence of adaptive control schemes using least-squares parameter
  estimates.
\newblock {\em IEEE Transactions on Automatic Control}, 35(4):416--424, 1990.

\bibitem{Lai86}
T.L. Lai.
\newblock Asymptotically efficient adaptive control in stochastic regression
  models.
\newblock {\em Advances in Applied Math.}, 7(23):23--45, 1986.

\bibitem{Lai94}
T.L. Lai.
\newblock Asymptotic properties of nonlinear least squares estimates in
  stochastic regression models.
\newblock {\em Annals of Statistics}, 22(4):1917--1930, 1994.

\bibitem{LaiR81}
T.L. Lai and H.~Robbins.
\newblock Consistency and asymptotic efficiency of slope estimates in
  stochastic approximation schemes.
\newblock {\em Z. Wahrsch. verw. Gebiete}, 56:329--360, 1981.

\bibitem{LaiRW78}
T.L. Lai, H.~Robbins, and C.Z. Wei.
\newblock Strong consistency of least squares estimates in multiple regression.
\newblock {\em Proc.\ Nat.\ Acad.\ Sci.\ USA}, 75(7):3034--3036, 1978.

\bibitem{LaiRW79}
T.L. Lai, H.~Robbins, and C.Z. Wei.
\newblock Strong consistency of least squares estimates in multiple regression
  {II}.
\newblock {\em Journal of Multivariate Analysis}, 9:343--361, 1979.

\bibitem{LaiW82}
T.L. Lai and C.Z. Wei.
\newblock Least squares estimates in stochastic regression models with
  applications to identification and control of dynamic systems.
\newblock {\em Annals of Statistics}, 10(1):154--166, 1982.

\bibitem{LaiW86}
T.L. Lai and C.Z. Wei.
\newblock On the concept of excitation in least squares identification and
  adaptive control.
\newblock {\em Stochastics}, 16:227--254, 1986.

\bibitem{LaiW87}
T.L. Lai and C.Z. Wei.
\newblock Asymptotically efficient self-tuning regulators.
\newblock {\em SIAM J. Control and Optimization}, 25(2):466--481, 1987.

\bibitem{LearyBK2004}
S.~Leary, A.~Bhaskar, and A.J. Keane.
\newblock A derivative based surrogate model for approximating and optimizing
  the output of an expensive computer simulation.
\newblock {\em Journal of Global Optimization}, 30:39--58, 2004.

\bibitem{LiangZ95}
K.-Y. Liang and S.L. Zeger.
\newblock Inference based on estimating functions in the presence of nuisance
  parameters.
\newblock {\em Statistical Science}, 10(2):158--173, 1995.

\bibitem{Ljung87}
L.~Ljung.
\newblock {\em System Identification, Theory for the User}.
\newblock Prentice-Hall, Englewood Cliffs, 1987.

\bibitem{Mardia90}
K.V. Mardia.
\newblock Maximum likelihood estimation for spatial models.
\newblock In D.A. Griffith, editor, {\em Spatial Statistics: Past, Present and
  Future}, pages 203--253. Michigan Document Services, Ann Arbor, Michigan,
  1990.

\bibitem{MardiaKGL96}
K.V. Mardia, J.T. Kent, C.R. Goodall, and J.A. Little.
\newblock Kriging and splines with derivative information.
\newblock {\em Biometrika}, 83(1):207--221, 1996.
\newblock (correction in Biometrika (1998), 85(2), p.\ 205).

\bibitem{Matheron63}
G.~Matheron.
\newblock Principles of geostatistics.
\newblock {\em Economic Geology}, 58:1246--1266, 1963.

\bibitem{Mitchell74}
T.J. Mitchell.
\newblock An algorithm for the construction of ``{$D$}-optimal'' experimental
  designs.
\newblock {\em Technometrics}, 16:203--210, 1974.

\bibitem{Mockus89}
J.~Mockus.
\newblock {\em {B}ayesian Approach to Global Optimization, Theory and
  Applications}.
\newblock Kluwer, Dordrecht, 1989.

\bibitem{MockusTZ78}
J.~Mockus, V.~Tiesis, and A.~Zilinskas.
\newblock The application of {B}ayesian methods for seeking the extremum.
\newblock In L.C.W. Dixon and G.P. Szego, editors, {\em Towards Global
  Optimisation 2}, pages 117--129. North Holland, Amsterdam, 1978.

\bibitem{MolchanovZ01}
I.~Molchanov and S.~Zuyev.
\newblock Variational calculus in the space of measures and optimal design.
\newblock In A.~Atkinson, B.~Bogacka, and A.~Zhigljavsky, editors, {\em Optimum
  Design 2000}, chapter~8, pages 79--90. Kluwer, Dordrecht, 2001.

\bibitem{MolchanovZ02}
I.~Molchanov and S.~Zuyev.
\newblock Steepest descent algorithm in a space of measures.
\newblock {\em Statistics and Computing}, 12:115--123, 2002.

\bibitem{MorrisM95}
M.D. Morris and T.J. Mitchell.
\newblock Exploratory designs for computational experiments.
\newblock {\em Journal of Statistical Planning and Inference}, 43:381--402,
  1995.

\bibitem{MorrisMY93}
M.D. Morris, T.J. Mitchell, and D.~Ylvisaker.
\newblock Bayesian design and analysis of computer experiments: use of
  derivatives in surface prediction.
\newblock {\em Technometrics}, 35(3):243--255, 1993.

\bibitem{Muller84}
H.-G. {M\"{u}ller}.
\newblock Optimal designs for nonparametric kernel regression.
\newblock {\em Statistics \& Probability Letters}, 2:285--290, 1984.

\bibitem{Muller2001}
W.G. {M\"{u}ller}.
\newblock {\em Collecting Spatial Data. Optimum Design of Experiments for
  Random Fields (2nd revised edition)}.
\newblock Physica-Verlag, Heidelberg, 2001.

\bibitem{MullerP2003}
W.G. {M\"{u}ller} and A.~P\'azman.
\newblock Measures for designs in experiments with correlated errors.
\newblock {\em Biometrika}, 90(2):423--434, 2003.

\bibitem{Nadaraya64}
E.A. Nadaraya.
\newblock On estimating regression.
\newblock {\em Theory of Probability and its Applications}, 9:141--142, 1964.

\bibitem{Nather75}
W.~N{\"a}ther.
\newblock The choice of estimators and experimental designs in a linear
  regression model according to a joint criterion of optimality.
\newblock {\em Math. Operationsforsch. Statist.}, 6:677--686, 1975.

\bibitem{NesterovN94}
Y.~Nesterov and A.~Nemirovskii.
\newblock {\em Interior-Point Polynomial Algorithms in Convex Programming}.
\newblock SIAM, Philadelphia, 1994.

\bibitem{NesicL2002}
D.~Ne\v{s}i\'c and D.S. Laila.
\newblock A note on input-to-state stabilization for nonlinear sampled-data
  systems.
\newblock {\em IEEE Transactions on Automatic Control}, 47(7):1153--1158, 2002.

\bibitem{NesicT2001}
D.~Ne\v{s}i\'c and A.R. Teel.
\newblock Sampled-data control of nonlinear systems: an overview of recent
  results.
\newblock In {\em Perspectives in robust control}, pages 221--239. Lecture
  Notes in Control and Inform. Sci., 268, Springer, New York, 2001.

\bibitem{Ninness2000}
B.~Ninness.
\newblock Strong laws of large numbers under weak assumptions with application.
\newblock {\em IEEE Transactions on Automatic Control}, 45(11):2117--2122,
  2000.

\bibitem{NinnessG95}
B.~Ninness and G.~Goodwin.
\newblock Estimation of model quality.
\newblock {\em Automatica}, 31(12):1771--1797, 1995.

\bibitem{Pazman86}
A.~P\'azman.
\newblock {\em Foundations of Optimum Experimental Design}.
\newblock Reidel (Kluwer group), Dordrecht (co-pub. VEDA, Bratislava), 1986.

\bibitem{Pazman93}
A.~P\'azman.
\newblock {\em Nonlinear Statistical Models}.
\newblock Kluwer, Dordrecht, 1993.

\bibitem{PazmanM2001}
A.~P\'azman and W.G. M{\"u}ller.
\newblock Optimum design of experiments subject to correlated errors.
\newblock {\em Statistics \& Probability Letters}, 52:29--34, 2001.

\bibitem{PazmanPa92}
A.~P\'azman and L.~Pronzato.
\newblock Nonlinear experimental design based on the distribution of
  estimators.
\newblock {\em Journal of Statistical Planning and Inference}, 33:385--402,
  1992.

\bibitem{PazmanPa96}
A.~P\'azman and L.~Pronzato.
\newblock A {D}irac function method for densities of nonlinear statistics and
  for marginal densities in nonlinear regression.
\newblock {\em Statistics \& Probability Letters}, 26:159--167, 1996.

\bibitem{PPmoda72004}
A.~P{\'a}zman and L.~Pronzato.
\newblock Simultaneous choice of design and estimator in nonlinear regression
  with parameterized variance.
\newblock In A.~{Di Bucchianico}, H.~L{\"a}uter, and H.P. Wynn, editors, {\em
  mODa'7 -- Advances in Model--Oriented Design and Analysis, Proceedings of the
  7th Int.\ Workshop, Heeze (Netherlands)}, pages 117--124, Heidelberg, June
  2004. Physica Verlag.

\bibitem{PPa05_SPL}
A.~P\'azman and L.~Pronzato.
\newblock On the irregular behavior of {LS} estimators for asymptotically
  singular designs.
\newblock {\em Statistics \& Probability Letters}, 76:1089--1096, 2006.

\bibitem{Pilz83}
J.~Pilz.
\newblock {\em {B}ayesian Estimation and Experimental Design in Linear
  Regression Models}, volume~55.
\newblock Teubner-Texte zur Mathematik, Leipzig, 1983.
\newblock (also {W}iley, {N}ew {Y}ork, 1991).

\bibitem{PoggiP2000}
J.-M. Poggi and B.~Portier.
\newblock Nonlinear adaptive tracking using kernel estimators: estimation and
  test for linearity.
\newblock {\em SIAM J.\ Control Optim.}, 39(3):707--727, 2000.

\bibitem{PortierO2000}
B.~Portier and A.~Oulidi.
\newblock Nonparametric estimation and adaptive control of functional
  autoregressive models.
\newblock {\em SIAM J.\ Control Optim.}, 39(2):411--432, 2000.

\bibitem{Pronzato00a}
L.~Pronzato.
\newblock Adaptive optimisation and {$D$}-optimum experimental design.
\newblock {\em Annals of Statistics}, 28(6):1743--1761, 2000.

\bibitem{Pa05}
L.~Pronzato.
\newblock On the sequential construction of optimum bounded designs.
\newblock {\em Journal of Statistical Planning and Inference}, 136:2783--2804,
  2006.

\bibitem{PKWa96}
L.~Pronzato, C.~Kulcs\'ar, and E.~Walter.
\newblock An actively adaptive control policy for linear models.
\newblock {\em IEEE Transactions on Automatic Control}, 41(6):855--858, 1996.

\bibitem{PPa94}
L.~Pronzato and A.~P\'azman.
\newblock Second-order approximation of the entropy in nonlinear least-squares
  estimation.
\newblock {\em Kybernetika}, 30(2):187--198, 1994.
\newblock {\em Erratum} 32(1):104, 1996.

\bibitem{PP01a}
L.~Pronzato and A.\ P\'azman.
\newblock Using densities of estimators to compare pharmacokinetic experiments.
\newblock {\em Computers in Biology and Medicine}, 31(3):179--195, 2001.

\bibitem{PP-eusipco-04}
L.~Pronzato and A.~P\'azman.
\newblock Recursively re-weighted least-squares estimation in regression models
  with parameterized variance.
\newblock In {\em Proc.\ EUSIPCO'2004, Vienna, Austria}, pages 621--624,
  September 2004.

\bibitem{PTa03}
L.~Pronzato and E.~Thierry.
\newblock Sequential experimental design and response optimisation.
\newblock {\em Statistical Methods and Applications}, 11(3):277--292, 2003.

\bibitem{PWa85}
L.~Pronzato and E.~Walter.
\newblock Robust experiment design via stochastic approximation.
\newblock {\em Mathematical Biosciences}, 75:103--120, 1985.

\bibitem{PWa88}
L.~Pronzato and E.~Walter.
\newblock Robust experiment design via maximin optimization.
\newblock {\em Mathematical Biosciences}, 89:161--176, 1988.

\bibitem{PWa93}
L.~Pronzato and E.~Walter.
\newblock Experimental design for estimating the optimum point in a response
  surface.
\newblock {\em Acta Applicandae Mathematicae}, 33:45--68, 1993.

\bibitem{PWa94_2}
L.~Pronzato and E.~Walter.
\newblock Minimum-volume ellipsoids containing compact sets: application to
  parameter bounding.
\newblock {\em Automatica}, 30(11):1731--1739, 1994.

\bibitem{PWZl2000}
L.~Pronzato, H.P. Wynn, and A.A. Zhigljavsky.
\newblock {\em Dynamical Search}.
\newblock Chapman \& Hall/CRC, Boca Raton, 2000.

\bibitem{PWZa01-AAM}
L.~Pronzato, H.P. Wynn, and A.A. Zhigljavsky.
\newblock Renormalised steepest descent in {H}ilbert space converges to a
  two-point attractor.
\newblock {\em Acta Applicandae Mathematicae}, 67:1--18, 2001.

\bibitem{PWZa05_MP}
L.~Pronzato, H.P. Wynn, and A.A. Zhigljavsky.
\newblock Asymptotic behaviour of a family of gradient algorithms in
  $\mathbb{R}^d$ and {H}ilbert spaces.
\newblock {\em Mathematical Programming}, A107:409--438, 2006.

\bibitem{Pukelsheim93}
F.~Pukelsheim.
\newblock {\em Optimal Experimental Design}.
\newblock Wiley, New York, 1993.

\bibitem{PukelsheimR92}
F.~Pukelsheim and S.~Reider.
\newblock Efficient rounding of approximate designs.
\newblock {\em Biometrika}, 79(4):763--770, 1992.

\bibitem{Rafajlowicz83}
E.~Rafaj{\l}owicz.
\newblock Optimal experiment design for identification of linear distributed
  parameter systems: Frequency domain approach.
\newblock {\em IEEE Transactions on Automatic Control}, 28(7):806--808, 1983.

\bibitem{Rafajlowicz86b}
E.~Rafaj{\l}owicz.
\newblock Optimum choice of moving sensor trajectories for distributed
  parameter system identification.
\newblock {\em International Journal of Control}, 43(5):1441--1451, 1986.

\bibitem{RaynaudPW01a}
H.-F. Raynaud, L.~Pronzato, and E.~Walter.
\newblock Robust identification and control based on ellipsoidal parametric
  uncertainty descriptions.
\newblock {\em European J. of Control}, 6(3):245--257, 2000.

\bibitem{RobertazziS89}
T.G. Robertazzi and S.C. Schwartz.
\newblock An accelerated sequential algorithm for producing {$D$}-optimal
  designs.
\newblock {\em SIAM J. Sci. Stat. Comput.}, 10(2):341--358, 1989.

\bibitem{RojasWGF07}
C.R. Rojas, J.S. Welsh, G.C. Goodwin, and A.~Feuer.
\newblock Robust optimal experiment design for system identification.
\newblock {\em Automatica}, 43:993--1008, 2007.

\bibitem{SacksS88}
J.~Sacks and S.~Schiller.
\newblock Spatial designs.
\newblock In S.S. Gupta and J.O. Berger, editors, {\em Statistical Decision
  Theory and Related Topics IV}, volume~2, pages 385--399. Springer,
  Heidelberg, 1988.

\bibitem{SacksWMW89}
J.~Sacks, W.J. Welch, T.J. Mitchell, and H.P. Wynn.
\newblock Design and analysis of computer experiments.
\newblock {\em Statistical Science}, 4(4):409--435, 1989.

\bibitem{SantnerWN2003}
T.~Santner, B.J. Williams, and W.I. Notz.
\newblock {\em The Design and Analysis of Computer Experiments}.
\newblock Springer, Heidelberg, 2003.

\bibitem{Schaback2003}
R.~Schaback.
\newblock Mathematical results concerning kernel techniques.
\newblock In {\em Prep.\ 13th IFAC Symposium on System Identification,
  Rotterdam}, pages 1814--1819, August 2003.

\bibitem{Scheffe61}
H.~Scheff\'e.
\newblock Simultaneous interval estimates of linear fucntions of parameters.
\newblock {\em Bull.\ Inst.\ Internat.\ Statist.}, 38:245--253, 1961.

\bibitem{SchonlauWJ98}
M.~Schonlau, W.J. Welch, and D.R. Jones.
\newblock Global versus local search in constrained optimization of computer
  models.
\newblock In {\em New Developments and Applications in Experimental Design,
  Lecture Notes --- Monograph Series, vol.\ 34}, pages 11--25. IMS, Hayward,
  1998.

\bibitem{Schwabe87b}
R.~Schwabe.
\newblock On adaptive chemical balance weighting designs.
\newblock {\em Journal of Statistical Planning and Inference}, 17:209--216,
  1987.

\bibitem{ShewryW87}
M.C. Shewry and H.P. Wynn.
\newblock Maximum entropy sampling.
\newblock {\em Applied Statistics}, 14:165--170, 1987.

\bibitem{Shiryaev96}
A.N. Shiryaev.
\newblock {\em Probability}.
\newblock Springer, Berlin, 1996.

\bibitem{Sibson72}
R.~Sibson.
\newblock Discussion on a paper by {H}.{P}. {W}ynn.
\newblock {\em Journal of Royal Statistical Society}, B34:181--183, 1972.

\bibitem{SilvermanT80}
B.W. Silverman and D.M. Titterington.
\newblock Minimum covering ellipses.
\newblock {\em SIAM Journal Sci. Stat. Comput.}, 1(4):401--409, 1980.

\bibitem{Silvey80}
S.D. Silvey.
\newblock {\em Optimal Design}.
\newblock Chapman \& Hall, London, 1980.

\bibitem{SilveyTT78}
S.D. Silvey, D.M. Titterington, and B.~Torsney.
\newblock An algorithm for optimal designs on a finite design space.
\newblock {\em Commun. Statist.-Theor. Meth.}, A7(14):1379--1389, 1978.

\bibitem{SoderstromS81}
T.~S{\"o}derstr{\"o}m and P.~Stoica.
\newblock Comparison of some instrumental variable methods---consistency and
  accuracy aspects.
\newblock {\em Automatica}, 17(1):101--115, 1981.

\bibitem{SoderstromS83}
T.~S{\"o}derstr{\"o}m and P.~Stoica.
\newblock {\em Instrumental Variable Methods for System Identification}.
\newblock Springer, New York, 1983.

\bibitem{SoderstromS89}
T.~S{\"o}derstr{\"o}m and P.~Stoica.
\newblock {\em System Identification}.
\newblock Prentice Hall, New York, 1989.

\bibitem{Spokoinyi92}
V.G. Spokoinyi.
\newblock On asymptotically optimal sequential experimental design.
\newblock {\em Advances in Soviet Mathematics}, 12:135--150, 1992.

\bibitem{Stein99}
M.L. Stein.
\newblock {\em Interpolation of Spatial Data. Some Theory for Kriging}.
\newblock Springer, Heidelberg, 1999.

\bibitem{Sternby77}
J.~Sternby.
\newblock On consistency for the method of least squares using martingale
  theory.
\newblock {\em IEEE Transactions on Automatic Control}, 22(3):346--352, 1977.

\bibitem{Titterington75}
D.M. Titterington.
\newblock Optimal design: some geometrical espects of ${D}$-optimality.
\newblock {\em Biometrika}, 62(2):313--320, 1975.

\bibitem{Titterington76}
D.M. Titterington.
\newblock Algorithms for computing {$D$}-optimal designs on a finite design
  space.
\newblock In {\em Proc. of the 1976 Conference on Information Science and
  Systems}, pages 213--216, Baltimore, 1976. Dept. of Electronic Engineering,
  John Hopkins University.

\bibitem{Torsney83}
B.~Torsney.
\newblock A moment inequality and monotonicity of an algorithm.
\newblock In K.O. Kortanek and A.V. Fiacco, editors, {\em Proc. Int. Symp. on
  Semi-infinite Programming and Applications}, pages 249--260, Heidelberg,
  1983. Springer.

\bibitem{TseB-S73}
E.~Tse and Y.~Bar-{S}halom.
\newblock An actively adaptive control for linear systems with random
  parameters via the dual control approach.
\newblock {\em IEEE Transactions on Automatic Control}, 18(2):109--117, 1973.

\bibitem{TseB-SM73}
E.~Tse, Y.~Bar-{S}halom, and L.~{Meier~III}.
\newblock Wide-sense adaptive dual control for nonlinear stochastic systems.
\newblock {\em IEEE Transactions on Automatic Control}, 18(2):98--108, 1973.

\bibitem{Ucinski2005}
D.~Uci{\'n}ski.
\newblock {\em Optimal Measurement Methods for Distributed Parameter System
  Identification}.
\newblock CRC Press, Boca Raton, 2005.

\bibitem{vanderVaart96}
A.W. van~der {V}aart.
\newblock Maximum likelihood estimation under a spatial sampling scheme.
\newblock {\em Annals of Statistics}, 24(5):2049--2057, 1996.

\bibitem{VandenbergheBW98}
L.~Vandenberghe, S.~Boyd, and S.-P. Wu.
\newblock Determinant maximisation with linear matrix inequality constraints.
\newblock {\em SIAM Journal\ on Matrix Analysis and Applications},
  19(2):499--533, 1998.

\bibitem{Vapnik98}
V.N. Vapnik.
\newblock {\em Statistical Learning Theory}.
\newblock Wiley, New York, 1998.

\bibitem{Vapnik2000}
V.N. Vapnik.
\newblock {\em The Nature of Statistical Learning Theory}.
\newblock Springer, New York, 2000.
\newblock 2nd Edition.

\bibitem{Vazquez05}
E.~Vazquez.
\newblock Mod\'elisation comportementale de syst\`emes nonlin\'eaires
  multivariables par m\'ethodes \`a noyaux et applications.
\newblock Ph.D. Thesis, Universit\'e Paris XI, Orsay, France, May 2005.

\bibitem{VazquezW03}
E.~Vazquez and E.~Walter.
\newblock Multi-output support vector regression.
\newblock In {\em Prep.\ 13th IFAC Symposium on System Identification,
  Rotterdam}, pages 1820--1825, August 2003.

\bibitem{VazquezWF2005}
E.~Vazquez, E.~Walter, and G.~Fleury.
\newblock Intrinsic {K}riging and prior information.
\newblock {\em Applied Stochastic Models in Business and Industry},
  21:215--226, 2005.

\bibitem{WPl97}
E.~Walter and L.~Pronzato.
\newblock {\em Identification of Parametric Models from Experimental Data}.
\newblock Springer, Heidelberg, 1997.

\bibitem{Watson64}
G.S. Watson.
\newblock Smooth regression analysis.
\newblock {\em Sankhya, Series A}, 26:359--372, 1964.

\bibitem{Wright97}
S.~Wright, editor.
\newblock {\em Primal-Dual Interior-Point Methods}.
\newblock SIAM, Philadelphia, 1997.

\bibitem{Wu85}
C.F.J. Wu.
\newblock Asymptotic inference from sequential design in a nonlinear situation.
\newblock {\em Biometrika}, 72(3):553--558, 1985.

\bibitem{Wynn70}
H.P. Wynn.
\newblock The sequential generation of {$D$}-optimum experimental designs.
\newblock {\em Annals of Math. Stat.}, 41:1655--1664, 1970.

\bibitem{Wynn2004}
H.P. Wynn.
\newblock Maximum entropy sampling and general equivalence theory.
\newblock In A.~{Di Bucchianico}, H.~L{\"a}uter, and H.P. Wynn, editors, {\em
  mODa'7 -- Advances in Model--Oriented Design and Analysis, Proceedings of the
  7th Int.\ Workshop, Heeze (Netherlands)}, pages 211--218. Physica Verlag,
  Heidelberg, June 2004.

\bibitem{Ye97}
Y.~Ye.
\newblock {\em Interior-Point Algorithms: Theory and Analysis}.
\newblock Wiley, Chichester, 1997.

\bibitem{Ying93}
Z.~Ying.
\newblock Maximum likelihood estimation of parameters under a spatial sampling
  scheme.
\newblock {\em Annals of Statistics}, 21:1567--1590, 1993.

\bibitem{Zarrop79}
M.B. Zarrop.
\newblock {\em Optimal Experiment Design for Dynamic System Identification}.
\newblock Springer, Heidelberg, 1979.

\bibitem{ZhuZ2006}
Z.~Zhu and H.~Zhang.
\newblock Spatial sampling design under the infill asymptotic framework.
\newblock {\em Environmetrics}, 17(4):323--337, 2006.

\end{thebibliography}



\end{document}